\documentclass[a4paper]{article}
\usepackage{lgrind}
\usepackage{amssymb}
\usepackage{mathrsfs}
\usepackage{latexcad}\usepackage{float}
\usepackage{amsmath}
\usepackage{amsfonts,amssymb}
\usepackage{graphicx,color}
\usepackage{psfig}
\usepackage[toc,page,title,titletoc,header]{appendix}
\catcode`\@=11 \@addtoreset{equation}{section}
\def\theequation{\thesection.\arabic{equation}}
\catcode`\@=12

\newtheorem{Theorem}{Theorem}[section]

\newtheorem{Lemma}{Lemma}[section]

\newtheorem{Remark}{Remark}[section]

\newtheorem{Proposition}{Proposition}[section]

\def\2{{I \hskip -1.0mm I}}
\def\3{{I \hskip -1.0mm I\hskip -1.0mm I}}
\def\4{{I \hskip -0.9mm V}}
\def\6{{V \hskip -1.35mm I}}

\def\ss{\scriptstyle}

\title{A new geometric flow with rotational invariance}
\author{De-Xing Kong  $\quad$ and $\quad$
Qiang Ru
\\
Department of Mathematics, Zhejiang University\\ Hangzhou 310027,
China}
\date{ }
\begin{document}
\maketitle
\begin{abstract}
In this paper we introduce a new geometric flow with rotational
invariance and prove that, under this kind of flow, an arbitrary
smooth closed contractible hypersurface in the Euclidean space
$\mathbb{R}^{n+1}\;(n\ge 1)$ converges to $\mathbb{S}^n$ in the
$C^{\infty}$-topology as $t$ goes to the infinity. This result
covers the well-known theorem of Gage and Hamilton in \cite{g-h} for
the curvature flow of plane curves and the famous result of Huisken
in \cite{hui} on the flow by mean curvature of convex surfaces,
respectively.

\vskip 6mm


\vskip 3mm

\end{abstract}

\newpage
\baselineskip=7mm

\section{Introduction}
Since the last quarter of twentieth century, using partial
differential equations to formulate and solve geometric problems has
become a trend and a dominating force. A new area called geometric
analysis was born. When looking back at the history of geometric
analysis, one could see numerous success stories of utilizing
differential equations to tackle important problems in geometry,
topology and physics. Typical and important examples would include
Yau's solution to the Calabi conjecture using the complex
Monge-Ampere equation (see Yau \cite{yau1}), Schoen's solution of
the Yamabe conjecture (see Schoen \cite{Schoen}), Schoen-Yau's proof
of the positive mass conjecture (see Schoen-Yau \cite{S-Y}),
Donaldson's work on 4-dimensional smooth manifolds using the
Yang-Mills equation (see Donaldson \cite{Dona}), and recently,
Perelman's solution to the century-old Poincar$\acute{e}$ conjecture
using Hamilton's beautiful theory on the Ricci flow, which is just a
nonlinear version of the classical heat equation (see
\cite{p1}-\cite{p3}). However, despite all these success, the
equations studied and utilized in geometry so far are almost
exclusively of elliptic or parabolic type. With few exceptions,
hyperbolic equations have not yet found their way into the study of
geometric or topological problems. More recently, Kong et al
introduced the hyperbolic geometric flow which is a fresh start of
an attempt to introduce hyperbolic partial differential equations
into the realm of geometry (see \cite{k} or \cite{kl}). The kind of
flow is a very natural tool to understand the wave character of
metrics, the wave phenomenon of curvatures, the evolution of
manifolds and their structures (see \cite{dkl},
\cite{klw}-\cite{klx}).

In this paper, we introduce a new geometric flow with rotational
invariance. This flow is described by, formally a system of
parabolic partial differential equations, essentially a coupled
system of hyperbolic-parabolic partial differential equations with
rotational invariance. More precisely, let $\mathscr{S}_t$ be a
family of hypersurfaces in the $(n+1)$-dimensional Euclidean space
$\mathbb{R}^{n+1}$ with coordinates $(x_1,\cdots, x_{n+1})$, without
loss of generality, we may assume that the family of hypersurfaces
$\mathscr{S}_t$ is given by
\begin{equation}\label{1.1}x=x(t,\theta_1,\cdots,\theta_n),\end{equation}
where $x=(x_1,\cdots,x_{n+1})^T$ is a vector-valued smooth function
of $t$ and $\theta=(\theta_1,\cdots,\theta_n)$, the new flow
considered here is given by the following evolution equation
\begin{equation}\label{1.2}
\frac{\partial x}{\partial t}+\sum^n_{i=1}\frac{\partial
(f_{i}(|x|)x) }{\partial\theta_i}=\frac{x}{|x|}\Delta
|x|,\end{equation} where $f_i(\nu)\;(i=1,\cdots,n)$ are $n$ given
smooth functions, ${\displaystyle\Delta
=\sum^n_{i=1}\frac{\partial^2 }{\partial\theta_i^2}}$ is the
Laplacian operator, and $|\bullet|$ stands for the norm of the
vector $\bullet$ in $\mathbb{R}^{n+1}$. It is easy to verify that
the equation (\ref{1.2}) possesses the rotational invariance which
plays an important role in the present paper.

We are interested in the deformation of a smooth closed
contractible hypersurface $x=x_0(\theta_1,\cdots,\theta_n)$ under
the flow (\ref{1.2}), that is, we consider how the hypersurface
$x_0$ is smoothly deformed, say, embedded into a smooth family of
hypersurfaces depending on a time parameter. This can be reduced
to solve the Cauchy problem for (\ref{1.2}) with the initial data
\begin{equation}\label{1.3}
t=0:\;\;x=x_0(\theta_1,\cdots,\theta_n).\end{equation} Obviously,
in the present situation, $x_0=x_0(\theta_1,\cdots,\theta_n)$ is a
vector-valued periodic function, say, defined on $[0,1]^n$. In
Section 2, we shall prove
\begin{Theorem}
If $f\in C^{1}$,  $x_{0}\in L^{\infty}  and \;
|x_0(\theta_1,\cdots,\theta_n)|>0,$ then the Cauchy problem
(\ref{1.2}),  (\ref{1.3}) admits a unique global smooth solution on
$[0,\infty)\times \mathbb{R}^n$.
\end{Theorem}

In particular, the following theorem will be proved in Section 3.
\begin{Theorem}
Suppose that $f_i(\nu)$ are all constants, i.e., $f_i(\nu)\equiv c_i
\;(i=1,\cdots,n)$, suppose furthermore that
$x_0=x_0(\theta_1,\cdots,\theta_n)$ is a smooth vector-valued
periodic function with the period $[0,1]^n$, and satisfies
\begin{equation}\label{1.4}
|x_0(\theta_1,\cdots,\theta_n)|>0,\quad
\forall\;(\theta_1,\cdots,\theta_n)\in [0,1]^n.\end{equation} Then
the Cauchy problem (\ref{1.2})-(\ref{1.3})   has a global smooth
solution $x=x(t,\theta_1,\cdots,\theta_n)$, and the solution
satisfies
\begin{equation}\label{1.5}
|x(t,\theta_1,\cdots,\theta_n)|\longrightarrow \bar{r}_0 \triangleq
\int_{[0,1]^n}|x_0(\theta_1,\cdots,\theta_n)|d\theta_1\cdots
d\theta_n \qquad {\rm as} \;\;\; t\nearrow \infty.\end{equation}
Moreover, if the hypersurfaces undergo suitable homotheties, then
the normalized hypersurfaces converge to a sphere in the
$C^{\infty}$-topology as $t$ goes to the infinity.
\end{Theorem}

\begin{Remark} The geometric meaning of the result in Theorem 1.2 can be shown in the following
figure:
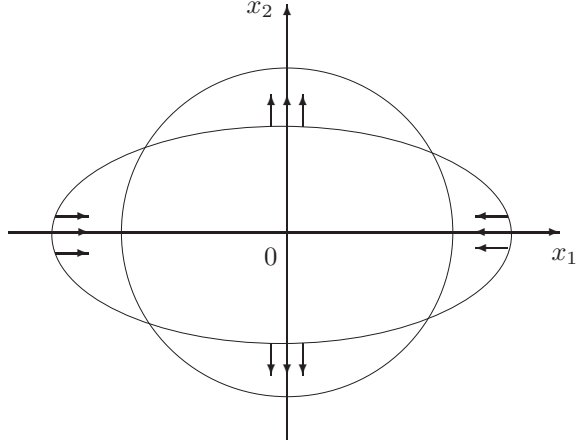
\begin{figure}[H]
    \begin{center}
\begin{picture}(198,174)
\thinlines \drawvector{-8.0}{76.0}{206.0}{1}{0}
\drawvector{96.0}{-2.0}{164.0}{0}{1}
\drawcenteredtext{90.0}{67.0}{0}
\drawcenteredtext{200.0}{67.0}{$x_1$}
\drawcenteredtext{86.0}{160.0}{$x_2$}
\drawcircle{96.0}{76.0}{124.0}{}
\drawellipse{94.0}{75.0}{172.0}{82.0}{}
\drawvector{9.6}{82.0}{12.0}{1}{0}
\drawvector{9.6}{68.0}{12.0}{1}{0}
\drawvector{8.0}{76.0}{14.0}{1}{0}
\drawvector{180.0}{76.0}{14.0}{-1}{0}
\drawvector{178.4}{82.0}{12.0}{-1}{0}
\drawvector{178.5}{70.0}{12.0}{-1}{0}
\drawvector{96.0}{116.0}{12.0}{0}{1}
\drawvector{102.0}{116.0}{12.0}{0}{1}
\drawvector{90.0}{116.0}{12.0}{0}{1}
\drawvector{96.0}{34.0}{12.0}{0}{-1}
\drawvector{90.0}{34.0}{12.0}{0}{-1}
\drawvector{102.0}{34.0}{12.0}{0}{-1}
\end{picture}
\caption{Deformation of an ellipse}
    \end{center}
\end{figure}
\end{Remark}

\begin{Remark} In the case $n=1$, Theorem 1.1 covers the well-known theorem of Gage and Hamilton in
\cite{g-h} for the curvature flow of plane curves; while for the
case of general $n$, Theorem 1.1 covers the famous result of Huisken
in \cite{hui} on the flow by mean curvature of convex surfaces into
spheres. In particular, we would like to point out that, in Theorem
1.1, we do NOT require the assumption that the hypersurface is
convex.
\end{Remark}

The paper is organized as follows. In Section 2 we prove the global
existence and uniqueness of smooth solutions for the Cauchy problem
(\ref{1.2}), (\ref{1.3}); Sections 3 is devoted to the proof of
Theorem 1.2; In Section 4, we state conclusions obtained in the
present paper and give some open problems. In Appendix, we
investigate the time-asymptotic behavior of global smooth solutions
for the equation (\ref{1.2}).

\section{Global existence and uniqueness of smooth solutions}

This section is devoted to the global existence and uniqueness of
smooth solution of the following equation
\begin{equation}\label{2.1}
\frac{\partial x}{\partial t}+\sum^m_{i=1}\frac{\partial
(f_i(|x|)x) }{\partial\theta_i}=\frac{x}{|x|}\Delta
|x|,\end{equation} where $x=(x_1,\cdots,x_n)^T$ is the unknown
vector-valued function, $f(\nu)=(f_1(\nu),\cdots, f_m(\nu))^T$ is
a given smooth vector-valed function, ${\displaystyle\Delta
=\sum^m_{i=1}\frac{\partial^2 }{\partial\theta_i^2}}$ is the
Laplacian operator, and $|\bullet|$ stands for the norm of the
vector $\bullet$ in $\mathbb{R}^{n}$.

Let
\begin{equation}\label{2.2-1}
x=rP,\quad r=|x|,\quad P=(p_1,\cdots,p_n)\in \mathbb{S}^{n-1}.
\end{equation}
Then it is easy to verify that the equation (\ref{2.1}) can be
rewritten as
$$\frac{\partial r}{\partial t}+\sum^m_{i=1}\frac{\partial (f_i(r)r)
}{\partial\theta_i}=\Delta r\eqno{(2.1a)}$$ and
$$\frac{\partial P}{\partial t}+\sum^m_{i=1}\frac{\partial P
}{\partial\theta_i}f_i(r)=0\eqno{(2.1b)}$$ for smooth solutions.

We now consider the Cauchy problem for the equation (\ref{2.1}),
equivalently, the system (2.1a)-(2.1b) with initial data
\begin{equation}\label{2.3}
t=0:\;\; r=r_{0}(\theta),\quad P=P_{0}(\theta),
\end{equation}
where $r_0(\theta)$ is a given scalar function of $\theta$, and
$P_0(\theta)$ is a given vector-valued function of $\theta$. In what
follows, we first investigate the local existence of smooth solution
of the above Cauchy problem.

As the standard way, let $K(t,\theta)$ be the fundamental solution
associated with the operator ${\displaystyle
\frac{\partial}{\partial t}-\triangle}$. That is to say,
\begin{equation}\label{2.4}
K(t,\theta)=(4\pi
t)^{-\frac{n}{2}}\exp\left\{-\frac{|\theta|^{2}}{4t}\right\}.\end{equation}
Then the solution $r=r(t,\theta)$ of the Cauchy problem
\begin{equation}\label{2.5}
\left\{\begin{array}{l}{\displaystyle \frac{\partial r}{\partial
t}+\sum^m_{i=1}\frac{\partial (f_i(r)r)
}{\partial\theta_i}=\Delta r,}\vspace{2mm}\\
t=0:\;\; r=r_{0}(\theta)\end{array}\right.\end{equation} has the
following integral representation
\begin{equation}\label{2.6}
r(t,\theta)=K(t,\theta)*r_{0}(\theta)+\sum^{m}_{j=1}\int^{t}_{0}K_{\theta_{j}}
(t-s,\theta)*(f_{j}(r(s,\theta))r(s,\theta))ds,
\end{equation}
where $*$ denotes the convolution with the space variables. We have
\begin{Lemma}\label{2.1}
Assume that
\begin{equation}\label{2.7}
f\in C^{1},\quad r_{0}\in L^{\infty},\end{equation} then there
exists a positive constant $T$ such that the Cauchy problem (2.5)
admits a unique smooth solution $=r(t,\theta)$ on the strip
\begin{equation}\label{2.8}
\Pi_{T}=\{(t,\theta)\,|\, t\in [0,T],\;\theta\in
\mathbb{R}^{m}\},\end{equation} where
\begin{equation}\label{2.9}
T=\min\left\{\left(\frac{M\pi^{\frac{1}{2}}}{2H}\right)^{2},\;\left(\frac{\pi^{\frac{1}{2}}}{4Hm}\right)^{2}\right\},
\end{equation}
in which
\begin{equation}\label{2.10}
M=\|r_{0}(\theta)\|_{L^{\infty}},\quad
H=\max_{i,j=1,\cdots,m}\left\{\sup_{|r|\leq(m+1)M}|g_{i}(r(t,\theta))|,\quad
\sup_{|r|\leq(m+1)M}\left|\frac{\partial}{\partial
\theta_{j}}g_{i}(r(t,\theta))\right|\right\},
\end{equation}
here $g_j=rf_j(r)\;\;(j=1,\cdots,m)$.
\end{Lemma}

\vskip 2mm

\noindent{\bf Proof.} Set
\begin{equation}\label{2.11}
G_{T}=\left\{r:\,[0,T]\times \mathbb{R}^{m}\rightarrow
L^{\infty}(\mathbb{R}^{m})\left|\|r(t,\cdot)\|_{L^{\infty} }\leq
(m+1)M\quad {\rm for} \;\; t\in [0,T]\right.\right\}
\end{equation}
and let $\mathcal{T}$ be the following integral operator
\begin{equation}\label{2.12}
\mathcal
{T}r(t,\theta)=K(t,\theta)*r_{0}(\theta)+\sum^{m}_{j=1}\int^{t}_{0}K_{\theta_{j}}(t-s,\theta)*(f_{j}(r(s,\theta))
r(s,\theta))ds.
\end{equation}
The solution $r=r(t,\theta)$ can be obtained as the
$L^{\infty}$-limit of the sequence $\{r^{k}\}$ defined by
\begin{equation}\label{2.13}r^{0}(t,\theta)=K(t,\theta)*r_{0}(\theta), \quad r^{k+1}=\mathcal
{T}r^{k}\quad (n=0,1\cdots).\end{equation}

To prove the above statement, we first claim that, for any
$t\in[0,T]$, it holds that
\begin{equation}\label{2.14}
\|r^{k}(t,\theta)\|_{L^{\infty}}\leq (m+1)M,\quad\forall\; k\in
\{0,1,2,\cdots\}.
\end{equation}
In what follows, we prove (\ref{2.14}) by the method of induction.

When $k=0$, we have
\begin{equation}\label{2.15}
\|r^{0}(t,\theta)\|_{L^{\infty}}=\|K(t,\theta)*r_{0}(\theta)\|_{L^{\infty}}.
\end{equation}
By Young's inequality, we obtain
\begin{equation}\label{2.16}
\|r^{0}(t,\theta)\|_{L^{\infty}} \leq
\|K(t,\theta)\|_{L^{1}}\|r_{0}(\theta)\|_{L^{\infty}} =
\|r_{0}(\theta)\|_{L^{\infty}}=M\leq (m+1)M.
\end{equation}

Now we assume that $\|r^{k}(t,\theta)\|_{L^{\infty}}\leq (m+1)M\;\;
(k\in \mathbb{N})$ holds. We next prove
\begin{equation}\label{2.17}
\|r^{k+1}(t,\theta)\|_{L^{\infty}}\leq (m+1)M.\end{equation} In
fact,
\begin{equation}\label{2.18}\begin{array}{lll}
{\displaystyle \|r^{k+1}(t,\theta)\|_{L^{\infty}}}&=&
{\displaystyle\|\mathcal
{T}r^{k}(t,\theta)\|_{L^{\infty}}}\vspace{2mm}\\
& \leq & {\displaystyle
\|K(t,x)*r_{0}(\theta)\|_{L^{\infty}}+\sum^{m}_{j=1}\int^{t}_{0}\|K_{x_{j}}(t-s,\theta)*(f_{j}(r^{k}
(s,\theta))r^{k}(s,\theta))\|_{L^{\infty}}ds}\vspace{2mm}\\
&\leq& {\displaystyle
M+\sum^{m}_{j=1}\int^{t}_{0}\|K_{\theta_{j}}(t-s,\theta)\|_{L^{1}}\|(f_{j}(r^{k}(s,\theta))r^{k}(s,\theta))
\|_{L^{\infty}}ds.}
\end{array}\end{equation}
Notice that
\begin{equation}\label{2.19}
\begin{array}{lll} \int_{\mathbb{R}^{m}}|
{\displaystyle K_{\theta_{j}}(t-s,\theta)|d\theta} & = &
{\displaystyle \int_{\mathbb{R}^{m}}[4\pi
(t-s)]^{-\frac{m}{2}}\frac{\theta_{j}}{2(t-s)}\exp
\left\{-\frac{|\theta|^{2}}{4(t-s)}\right\}d\theta}\vspace{2mm}\\
&=&{\displaystyle [4\pi
(t-s)]^{-\frac{m}{2}}\int_{\mathbb{R}^{m-1}}\left\{\exp
\left\{-\frac{\theta_{1}^{2}}{4(t-s)}\right\}+\cdots+\exp
\left\{-\frac{\theta_{j-1}^{2}}{4(t-s)}
\right\}+\right.}\vspace{2mm}\\
& & {\displaystyle
\left.\exp\left\{-\frac{\theta_{j+1}^{2}}{4(t-s)}\right\}+\cdots+
\exp\left\{-\frac{\theta_{m}^{2}}{4(t-s)}\right\}\right\}
d\theta_{1}\cdots d\theta_{j-1} d\theta_{j+1}\cdots
d\theta_{m}\times }\vspace{2mm}\\ & &
{\displaystyle\int_{\mathbb{R}}\frac{\theta_{j}}{2(t-s)}\exp
\left\{-\frac{\theta_{j}^{2}}{4(t-s)}\right\}d\theta_{j}}\vspace{2mm}\\
&=&{\displaystyle [4\pi
(t-s)]^{-\frac{m}{2}}\left\{[4(t-s)]^{\frac{1}{2}}\int_{\mathbb{R}}\exp
\left\{-\frac{\theta_{1}^{2}}{4(t-s)}\right\}d\left(-\frac{\theta_{1}}{[4(t-s)]^{\frac{1}{2}}}\right)\right\}^{m-1}
\times} \vspace{2mm}\\
& &{\displaystyle \int_{\mathbb{R}}\frac{\theta_{j}}{2(t-s)}\exp
\left\{-\frac{\theta_{j}^{2}}{4(t-s)}\right\}d\theta_{j}}\vspace{2mm}\\
&=&{\displaystyle [4\pi
(t-s)]^{-\frac{1}{2}}\int_{\mathbb{R}}\frac{\theta_{j}}{2(t-s)}\exp
\left\{-\frac{\theta_{j}^{2}}{4(t-s)}\right\}d\theta_{j}}\vspace{2mm}\\
&=&{\displaystyle [4\pi
(t-s)]^{-\frac{1}{2}}\left[-2\int_{0}^{\infty}\exp
\left\{-\frac{\theta_{j}^{2}}{4(t-s)}\right\}d\left(-\frac{\theta_{j}^{2}}{4(t-s)}\right)\right]}\vspace{2mm}\\
&=&{\displaystyle [\pi(t-s)]^{-\frac{1}{2}}.}
\end{array}\end{equation}
It follows from (\ref{2.18}) that
\begin{equation}\label{2.20}
\begin{array}{lll}{\displaystyle
\|r^{k+1}(t,\theta)\|_{L^{\infty}}} &\leq & {\displaystyle
M+\sum^{m}_{j=1}\int^{t}_{0}\|K_{\theta_{j}}(t-s,\theta)\|_{L^{1}}\|(f_{j}(r^{k}(s,\theta))r^{k}
(s,\theta))\|_{L^{\infty}}ds}\vspace{2mm}\\
&\leq& {\displaystyle
 M+\sum^{m}_{j=1}\pi^{-\frac{1}{2}}\int^{t}_{0}(t-s)^{-\frac{1}{2}}\|(f_{j}(r^{k}(s,\theta))r^{k}(s,\theta))
\|_{L^{\infty}}ds}\vspace{2mm}\\
&\leq& {\displaystyle  M+m\pi^{-\frac{1}{2}}H\int^{t}_{0}(t-s)^{-\frac{1}{2}}ds}\vspace{2mm}\\
&=&{\displaystyle M+2m\pi^{-\frac{1}{2}}Ht^{\frac{1}{2}}\leq
M+2m\pi^{-\frac{1}{2}}HT^{\frac{1}{2}}\leq (m+1)M.}
\end{array}\end{equation}
This is the desired estimate (\ref{2.14}). Thus, the proof of
(\ref{2.14}) is completed.

In what follows, we prove that $\{r^{k}(t,\theta)\}$ is uniformly
convergent in the strip $(0,T]\times \mathbb{R}^{m}$. To do so, it
suffices to show that $$
\sum^{\infty}_{k=1}\left[r^{k+1}(t,\theta)-r^{k}(t,\theta)\right]$$
is uniformly convergent in the strip $(0,T]\times \mathbb{R}^{m}$.

In fact, it holds that
\begin{equation}\label{2.21}
\begin{array}{lll}{\displaystyle
\|r^{k+1}-r^{k}\|_{L^{\infty}}} &\leq&{\displaystyle
\sum^{m}_{j=1}\int^{t}_{0}\|K_{\theta_{j}}(t-s,\theta)*\left[(f_{j}(r^{k}(s,\theta))r^{k}(s,\theta))-
(f_{j}(r^{k-1}(s,\theta))r^{k-1}(s,\theta))\right]\|_{L^{\infty}}ds}\vspace{2mm}\\
&\leq&{\displaystyle
\sum^{m}_{j=1}\int^{t}_{0}\|K_{\theta_{j}}(t-s,\theta)\|_{L^{1}}\|\left[(f_{j}(r^{k}(s,\theta))r^{k}(s,\theta))-
(f_{j}(r^{k-1}(s,\theta))r^{k-1}(s,\theta))\right]\|_{L^{\infty}}ds}\vspace{2mm}\\
&\leq&{\displaystyle
\sum^{m}_{j=1}\int^{t}_{0}\|K_{\theta_{j}}(t-s,\theta)\|_{L^{1}}\|g_{j}(r^{k}(s,\theta))-
g_{j}(r^{k-1}(s,\theta))\|_{L^{\infty}}ds}\vspace{2mm}\\
&\leq& {\displaystyle
\sum^{m}_{j=1}\int^{t}_{0}\|K_{\theta_{j}}(t-s,\theta)\|_{L^{1}}\|\nabla
g_{j}(\beta_{k})\|_{L^{\infty}}\|r^{k}(s,\theta)-r^{k-1}(s,\theta)\|_{L^{\infty}}ds}\vspace{2mm}\\
&\leq& {\displaystyle mH
\max\left\{|r^{k}(s,\theta)-r^{k-1}(s,\theta)|\right\}\int^{t}_{0}\|K_{\theta_{j}}(t-s,\theta)
\|_{L^{1}}ds}\vspace{2mm}\\
&\leq& {\displaystyle
2\pi^{-\frac{1}{2}}mHt^{\frac{1}{2}}\max\left\{|r^{k}(s,\theta)-r^{k-1}(s,\theta)|\right\}}
\vspace{2mm}\\
&\leq &{\displaystyle
2\pi^{-\frac{1}{2}}mHT^{\frac{1}{2}}\max\left\{|r^{k}(s,\theta)-r^{k-1}(s,\theta)|\right\}}
\vspace{2mm}\\
&\leq& {\displaystyle
\left(2\pi^{-\frac{1}{2}}mHT^{\frac{1}{2}}\right)^{2}\max\left\{|r^{k-1}(s,\theta)-r^{k-2}(s,\theta)|
\right\}}\vspace{2mm}\\
&\leq&\cdots\vspace{2mm}\\
& \leq & {\displaystyle
\left(2\pi^{-\frac{1}{2}}mHT^{\frac{1}{2}}\right)^{k}\max\left\{|r^{1}(s,\theta)-r^{0}(s,\theta)|
\right\},}
\end{array}\end{equation}
where
$$\beta_{k}\in\left[\min\{r^{k}(s,x),r^{k-1}(s,x)\},\max\{r^{k}(s,x),r^{k-1}(s,x)\}\right].$$
Noting
\begin{equation}\label{2.22}
\|r^{1}(s,\theta)-r^{0}(s,\theta)\|_{L^{\infty}}\leq
2\pi^{-\frac{1}{2}}mHT^{\frac{1}{2}},
\end{equation}
we obtain from (\ref{2.21}) that
\begin{equation}\label{2.23}
\|r^{k+1}-r^{k}\|_{L^{\infty}}\leq
\left(2\pi^{-\frac{1}{2}}mHT^{\frac{1}{2}}\right)^{k+1}.
\end{equation}
By (\ref{2.9}), we have
\begin{equation}\label{2.24}
\|r^{k+1}-r^{k}\|_{L^{\infty}}\leq \left(\frac12\right)^{k+1},
\end{equation}
which implies that
${\displaystyle\sum^{\infty}_{k=1}\left[r^{k+1}(t,\theta)-r^{k}(t,\theta)\right]}$
is uniformly convergent in the strip $(0,T]\times \mathbb{R}^{m}$.
Therefore, ${\displaystyle\lim_{k\rightarrow\infty}r^{k}(t,\theta)}$
gives the unique local solution of the Cauchy problem (2.5). Thus,
the proof Lemma 2.1 is completed.$\qquad\Box$

\begin{Lemma} Suppose that
$$f\in C^{1},\quad r_{0}\in L^{\infty}$$ and let
$M\stackrel{\triangle}{=}\|r_{0}\|_{L^{\infty}}$. Suppose
furthermore that $r(t,\theta)$ is the solution of Cauchy problem
(2.5) on the strip $\Pi_{T}$, then it holds that
\begin{equation}\label{2.25}
\|r(t,\theta)\|_{L^{\infty}(\Pi_{T}) }\leq M.
\end{equation}
\end{Lemma}

\vskip 2mm

\noindent{\bf Proof.} It follows from the proof of Lemma 2.1 that
\begin{equation}\label{2.26}
\|r(t,\theta)\|_{L^{\infty}(\Pi_{T}) }\leq (m+1)M
\stackrel{\triangle}{=}K.
\end{equation}
Introduce
\begin{equation}\label{2.27}
w(t,\theta)= r(t,\theta)-
M-\frac{K}{L^{2}}\left(|\theta|^{2}+CLe^{t}\right),
\end{equation}
where $C$ and $L$ are positive constants to be determined. By
(\ref{2.27}),
\begin{equation}\label{2.28}
r_{t}=w_{t}+\frac{CK}{L}e^{t},\quad \triangle r=\triangle
w+\frac{2Km}{L^{2}}.
\end{equation}
On the other hand,
\begin{equation}\label{2.29}
\sum^{m}_{j=1}(f_{j}(r)r)_{\theta_{j}}=\sum^{m}_{j=1}(g_{j}(r))_{\theta_{j}}
=\sum^{m}_{j=1}g^{\prime}_{j}(r)r_{\theta_{j}}
=\sum^{m}_{j=1}g^{\prime}_{j}(r)\left(w_{\theta_{j}}+\frac{2K}{L^{2}}\theta_{j}\right).
\end{equation}
Thus,
\begin{equation}\label{2.30}
w_{t}+\sum^{m}_{j=1}g^{\prime}_{j}(r)w_{\theta_{j}}+\frac{2K}{L^{2}}\sum^{m}_{j=1}g^{\prime}_{j}(r)\theta_{j}
+ \frac{CK}{L}e^{t}-\frac{2Km}{L^{2}}=\triangle w.
\end{equation}
Choose sufficiently large $C$ such that
\begin{equation}\label{2.31}
w(0,\theta)=r_{0}(\theta)-M-\frac{K}{L^{2}}(|\theta|^{2}+CL)<0,\quad
\forall\; \theta\in \mathbb{R}^m,
\end{equation}
and
\begin{equation}\label{2.32}\left\{\begin{array}{l} w(t,\pm L,\theta_{2},\cdots,
\theta_{m})=r(t,\pm L,\theta_{2},\cdots, \theta_{m})
M-\frac{K}{L^{2}}\left[(L^{2}+\theta_{2}^{2}+\cdots+\theta_{m}^{2})+CLe^{t}\right]<0,\vspace{2mm}\\
w(t,\theta_{1},\pm L,\cdots, \theta_{m})=r(t,\theta_{1},\pm
L,\cdots, \theta_{m})-
M-\frac{K}{L^{2}}\left[(\theta_{1}^{2}+L^{2}+\cdots+\theta_{m}^{2})+CLe^{t}\right]<0,\vspace{2mm}\\
\cdots\cdots\vspace{2mm}\\
w(t,\theta_{1},\cdots, \theta_{m-1},\pm L)=r(t,\theta_{1},\cdots,
\theta_{m-1},\pm L)
-M-\frac{K}{L^{2}}\left[(\theta_{1}^{2}+\cdots+\theta_{m-1}^{2}+|L|^{2})+CLe^{t}\right]<0
\end{array}\right.\end{equation}
for all $t\in[0,T]$.

In what follows, we prove that, for any $(t,\theta)\in
(0,T)\times(-L,L)^{m}$, it holds that
\begin{equation}\label{2.33}
w(t,\theta)<0.
\end{equation}
In fact, if (\ref{2.33}) is not true, then we can define $\bar{t}$
by
\begin{equation}\label{2.34}
\bar{t}=\inf_{t\in (0,T]}\{t\,|\,w(t,\theta)=0 \; \; {\rm for\;
some}\; \; \theta\in (-L,L)^{m}\}.
\end{equation}
It is easy to see that there exists a point, denoted by
$\bar{\theta}\in (-L,L)^{m}$, such that
\begin{equation}\label{2.35}
w(\bar{t},\bar{\theta})=0,\quad
w_{\theta_{1}}(\bar{t},\bar{\theta})=0,\quad \cdots,\quad
w_{\theta_{m}}(\bar{t},\bar{\theta})=0
\end{equation}
and
\begin{equation}\label{2.36}
w_{\theta_{i}\theta_{i}}(\bar{t},\bar{\theta})\leq 0,\quad\forall\;
i\in \{1,\cdots,m\}.
\end{equation}
By (\ref{2.35})-(\ref{2.36}), it follows from (\ref{2.30}) that
\begin{equation}\label{2.37}
w_{t}(\bar{t},\bar{\theta})+\frac{2K}{L^{2}}\sum^{m}_{j=1}g^{\prime}_{j}(r(\bar{t},\bar{\theta}))\bar{\theta_{j}}
+\frac{CK}{L}e^{\bar{t}}-\frac{2Km}{L^{2}}\leq 0.
\end{equation}
Noting
\begin{equation}\label{2.38}
\|g^{\prime}_{j}(\bullet)\|_{L^{\infty}} < \infty \quad {\rm
and}\quad (\bar{t},\theta_{j})\in (0,T]\times (-L,L),
\end{equation}
we can choose a sufficiently large $C$ such that
\begin{equation}\label{2.39}
\frac{2K}{L^{2}}\sum^{m}_{j=1}g^{\prime}_{j}(r(\bar{t},\bar{\theta}))\bar{\theta_{j}}
+\frac{CK}{L}e^{\bar{t}}-\frac{2Km}{L^{2}}>0.
\end{equation}
Combining (\ref{2.37}) and (\ref{2.39})
\begin{equation}\label{2.40}
w_{t}(\bar{t},\bar{\theta})<0.
\end{equation}
On the other hand, by the definition of $(\bar{t},\bar{\theta})$ it
holds that
\begin{equation}\label{2.41}
w_{t}(\bar{t},\bar{\theta})=\lim_{\triangle t\rightarrow
0}\frac{w(\bar{t},\bar{\theta})-w(\bar{t}-\triangle
t,\bar{\theta})}{\triangle t}\geq 0,
\end{equation}
which is a contradiction. This proves (\ref{2.33}).

Noting (\ref{2.27}) and (\ref{2.33}) and letting
$L\rightarrow\infty$ gives
\begin{equation}\label{2.42}
r(t,\theta)\leq M,\quad \forall\; (t,\theta)\in \Pi_{T}.
\end{equation}
Similarly, letting
\begin{equation}\label{2.43}
w(t,\theta)= r(t,\theta)+
M+\frac{K}{L^{2}}\left(|\theta|^{2}+CLe^{t}\right),\end{equation} we
can prove
\begin{equation}\label{2.44}
r(t,\theta)\geq -M,\quad \forall\; (t,\theta)\in \Pi_{T}.
\end{equation}
Combining (\ref{2.42}) and (\ref{2.44}) leads to
\begin{equation}\label{2.45}
\|r(t,\theta)\|_{L^{\infty}(\Pi_{T})}\leq M
\end{equation}
Thus, the proof of Lemma 2.2 is completed. $\qquad\Box$

\vskip 4mm

By Lemma 2.1 and Lemma 2.2, we have
\begin{Theorem}
If $f\in C^{1}$ and $r_{0}\in L^{\infty}$, then the Cauchy problem
(2.5) admits a unique global smooth solution on $[0,\infty)\times
\mathbb{R}^m$.
\end{Theorem}

Now we turn to consider the Cauchy problem (\ref{2.1}) (i.e.,
(2.1a)-(2.1b)), (\ref{2.3}). We have
\begin{Theorem}
Under the assumptions of Theorem 1.1, the Cauchy problem
(\ref{2.1}), (\ref{2.3}) admits a unique global smooth solution on
$[0,\infty)\times \mathbb{R}^m$.
\end{Theorem}

\vskip 2mm

\noindent{\bf Proof.}  Noting (1.4), by maximum principle we obtain
that, on the existence domain of smooth solution, it holds that
\begin{equation}\label{2.46}
r(t,\theta)> 0.
\end{equation}

On the one hand, we observe that, under the condition (\ref{2.46}),
the equation (\ref{2.1}) can be reduced to the system (2.1a)-(2.1b);
on the other hand, we notice that, once $r=r(t,\theta)$ is solved
from the Cauchy problem (\ref{2.5}), then the equation (2.1b)
becomes linear. Therefore, Theorem 2.2 follows Theorem 2.1 directly.
$\qquad\Box$

Obviously, Theorem 1.1 follows Theorem 2.2 directly.

\section{Time-asymptotic behavior of smooth solutions --- Proof of Theorem 1.2}

In this section, we first investigate the time-asymptotic behavior
of solution of the Cauchy problem (\ref{1.2}), (\ref{1.3}) in the
case that $f_{i}(\nu)\equiv c_{i}$, and then based on this, we prove
Theorem 1.2.

Notice that, in the present situation, the equation (2.1a) can be
rewritten as
\begin{equation*}
\frac{\partial r}{\partial t}+\sum^m_{i=1}c_{i}\frac{\partial r
}{\partial\theta_i}=\Delta r.
\end{equation*}
Let
\begin{equation*}
\tau=t,  \;\;\eta_{i}=\theta_{i}-c_{i}t,
\end{equation*}
then
\begin{equation*}
\frac{\partial  }{\partial\tau}=\frac{\partial}{\partial
t}\frac{\partial t}{\partial \tau}+\sum^m_{i=1}\frac{\partial
}{\partial\theta_i}\frac{\partial \theta_i}{\partial
\tau}=\frac{\partial}{\partial t}+\sum^m_{i=1}c_{i}\frac{\partial
}{\partial\theta_i}.
\end{equation*}
Therefore, without loss of generality, we turn to consider the
time-asymptotic behavior of solution of the Cauchy problem
\begin{equation}\label{3.1}
\left\{\begin{array}{lll} &\!\!\!\!\!\! u_{t}-\Delta u=0 \quad &
\textrm{in}\;\;(0,\infty)\times {\mathbb{R}}^n,\vspace{2mm} \\
&\!\!\!\!\!\!  u(x,0)= u_0 (x)  \quad &\textrm{in}\;\; \mathbb{R}^n,
\end{array}
\right.
\end{equation}
where $u=u(t,x)$ is the unknown function, the initial data $u_0(x)
\in C(T^{n})$ is a periodic function with period, say,
$\mathbb{T}^{n}=\{x=(x_{1},\cdots, x_{n})\mid -l_{i}\leq x_{i}\leq
l_{i}\}$.
\begin{Theorem}
If the initial data $u_0(x)$ satisfies Dirichlet conditions, then it
holds that
\begin{equation}\label{3.2}
u(t,x)\longrightarrow  \overline{u_0}\quad {\text as} \;\; t\nearrow
\infty,\end{equation} where $\overline{u_0}$ stands for the mean
value of $u_0(x)$, which is defined by
\begin{equation}\label{3.3}
\overline{u_0}=\frac{1}{{\rm
vol}\{\mathbb{T}^n\}}\int_{\mathbb{T}^n}u_0(x)dx.\end{equation}
\end{Theorem}
\vskip 2mm \noindent{\bf Proof.} Recall Dirichlet conditions: if the
periodic function $u_0(x)$ satisfies Dirichlet conditions, then
\begin{itemize}
  \item $u_0(x)$ has a finite number of extrema in any given interval;
\vskip -6mm
  \item  $u_0(x)$ has a finite number of discontinuities in any given
interval;
\vskip -6mm
  \item $u_0(x)$ is absolutely integrable over a period;
\vskip -6mm
  \item  $u_0(x)$ is bounded.
 \end{itemize}
\noindent By the theory of Fourier series, under Dirichlet
conditions, $u_0(x)$ is equal to the sum of its Fourier series at
each point where  $u_0(x)$ is continuous; moreover, the behavior of
the Fourier series at points of discontinuity is determined as well.
Therefore, it holds that
\begin{equation}\label{3.4}\begin{array}{lll}
u_0(x) &=&{\displaystyle
\sum^{\infty}_{m_{1}=0}\cdots\sum^{\infty}_{m_{n}=0}A^{0}_{m_{1}\cdots
m_{n}}\cos\frac{m_{1}\pi }{l_{1}}x_{1}\cdots \cos\frac{m_{n}\pi }{l_{n}}x_{n}+}\vspace{2mm}\\
& &{\displaystyle \sum^{n}_{i=1}\left\{\sum^{\infty}_{m_{1}=0}\cdots
\sum^{\infty}_{m_{i}=1}
\cdots\sum^{\infty}_{m_{n}=0}A^{1}_{m_{1}\cdots m_i\cdots m_{n}}\cos
\frac{m_{1}\pi }{l_{1}}x_{1}\cdots\sin\frac{m_{i}\pi }{l_{i}}x_{i}\cdots \cos\frac{m_{n}\pi }{l_{n}}x_{n}\right\}+}\vspace{2mm}\\
& &{\displaystyle
\sum^{n-1}_{i=1}\left\{\sum^{n}_{j=i+1}\left\{\sum^{\infty}_{m_{1}=0}\cdots
\sum^{\infty}_{ m_{i}=1} \cdots\sum^{\infty}_{ \stackrel{\ss
m_{j}=1}{i<j}}\cdots\sum^{\infty}_{m_{n}=0}A^{2}_{m_{1}\cdots
m_{i}\cdots
m_j\cdots m_{n}}\times \right.\right.}\vspace{2mm}\\
& &{\displaystyle \left.\left.\quad \cos \frac{m_{1}\pi
}{l_{1}}x_{1}\cdots\sin\frac{m_{i}\pi }{l_{i}}
x_{i}\cdots\sin\frac{m_{j}\pi }{l_{j}}x_{j}\cdots\cos \frac{m_{n}\pi }{l_{n}}x_{n}\right\}\right\}+\cdots +}\vspace{2mm}\\
& &{\displaystyle \sum^{\infty}_{m_{1}=1}\cdots
\sum^{\infty}_{m_{n}=1}A^{n}_{m_{1}\cdots m_{n}}\sin \frac{m_{1}\pi
}{l_{1}}x_{1}\cdots \sin \frac{m_{n}\pi }{l_{n}}x_{n}},
\end{array}\end{equation}
where $A^{j}_{m_{1}\cdots m_{n}}\; (j=0,1,2\cdots n)$ stand for the
Fourier coefficients which are given by
\begin{equation}\label{3.5}\left\{\begin{array}{lll}
A^{0}_{m_{1}\cdots m_{n}}&=&{\displaystyle
\frac{\lambda_{m_{1}\cdots m_{n}}}{l_{1}\cdots
l_{n}}\int^{l_{1}}_{-l_{1}}\cdots\int^{l_{n}}_{-l_{n}}u_0(x)\cos
\frac{m_{1}\pi }{l_{1}}x_{1}\cdots \cos \frac{m_{n}\pi
}{l_{n}}x_{n}dx_{1}\cdots dx_{n},}\vspace{2mm}\\
A^{1}_{m_{1}\cdots m_i\cdots m_{n}}&=&{\displaystyle
\frac{\lambda_{m_{1}\cdots m_{n}}}{l_{1}\cdots
l_{n}}\int^{l_{1}}_{-l_{1}}\cdots\int^{l_{i}}_{-l_{i}}\cdots\int^{l_{n}}_{-l_{n}}u_0(x)\times} \vspace{2mm}\\
&&{\displaystyle \cos \frac{m_{1}\pi
}{l_{1}}x_{1}\cdots\sin\frac{m_{i}\pi }{l_{i}}x_{i}\cdots \cos
\frac{m_{n}\pi }{l_{n}}x_{n}dx_{1}\cdots dx_i\cdots dx_{n},}
\vspace{2mm}\\
A^{2}_{m_{1}\cdots m_{i}\cdots m_j\cdots m_{n}}&=&{\displaystyle
\frac{\lambda_{m_{1}\cdots m_{n}}}{l_{1}\cdots
l_{n}}\int^{l_{1}}_{-l_{1}}\cdots\int^{l_{i}}_{-l_{i}}\cdots
\int^{l_{j}}_{-l_{j}} \cdots\int^{l_{n}}_{-l_{n}}u_0(x)\times }\vspace{2mm}\\
& &{\displaystyle \cos \frac{m_{1}\pi
}{l_{1}}x_{1}\cdots\sin\frac{m_{i}\pi
}{l_{i}}x_{i}\cdots\sin\frac{m_{j}\pi }{l_{j}}x_{j}\cdots\cos
\frac{m_{n}\pi }{l_{n}}x_{n}dx_{1}\cdots dx_{i}\cdots dx_{j}\cdots dx_{n},}\vspace{2mm}\\
\cdots\cdots & &  \vspace{2mm}\\
A^{n}_{m_{1}\cdots m_{n}}&=&{\displaystyle
\frac{\lambda_{m_{1}\cdots m_{n}}}{l_{1}\cdots
l_{n}}\int^{l_{1}}_{-l_{1}}\cdots\int^{l_{n}}_{-l_{n}}u_0(x)\sin
\frac{m_{1}\pi }{l_{1}}x_{1}\cdots \sin \frac{m_{n}\pi
}{l_{n}}x_{n}dx_{1}\cdots dx_{n}},
\end{array}\right.\end{equation}
in which
\begin{equation}\label{3.6}\left\{\begin{array}{l}
{\displaystyle \lambda_{m_{1}\cdots m_{n}}=\frac{1}{2^{n}},\quad
{\rm if}\;\;
m_{1}=m_{2}=\cdots = m_{n}=0,}\vspace{2mm} \\
{\displaystyle \lambda_{m_{1}\cdots m_{n}}=\frac{1}{2^{n-1}},\quad
{\rm if}\;\;
m_{i_{1}}=m_{i_{2}}=\cdots =m_{i_{n-1}}=0,\;\; m_{i_{n}}\neq 0,}\vspace{2mm} \\
\cdots\cdots\, \vspace{2mm} \\
{\displaystyle\lambda_{m_{1}\cdots m_{n}}=1,\quad {\rm if}\;\;
m_{1}=m_{2}=\cdots= m_{n}\neq 0.}
\end{array}\right.
\end{equation}
It is easy to see that
\begin{equation}\label{3.7}\begin{array}{lll}
u_0(x)&=&{\displaystyle\frac{1}{2^{n}l_{1}\cdots
l_{n}}\int^{l_{1}}_{-l_{1}}\cdots\int^{l_{n}}_{-l_{n}}u_0(x)dx_{1}\cdots
dx_{n}+\sum^{n}_{i=1}\left\{ \sum^{\infty}_{m_{i}=1} A^{0}_{0\cdots
m_{i} \cdots 0 }\cos \frac{m_{i}\pi
}{l_{i}}x_{i}\right\}+}\vspace{2mm}\\
&&{\displaystyle\sum^{n-1}_{i=1}\left\{\sum^{n}_{j=i+1}\left\{
\sum^{\infty}_{ m_{i}=1} \sum^{\infty}_{\stackrel{\ss m_{j}=1}
{i<j}}A^{0}_{0\cdots  m_{i}\cdots  m_{j} \cdots0 }\cos
\frac{m_{i}\pi }{l_{i}}x_{i}\cos\frac{m_{j}\pi
}{l_{j}}x_{j}\right\}\right\}+}\vspace{2mm}\\
&&{\displaystyle
\sum^{n-2}_{i=1}\left\{\sum^{n-1}_{j=i+1}\left\{\sum^{n}_{k=j+1}\left\{
\sum^{\infty}_{ m_{i}=1} \sum^{\infty}_{\stackrel{\ss m_{j}=1}{i<j}}
\sum^{\infty}_{\stackrel{\ss m_{k}=1} {i<j<k}}A^{0}_{0\cdots  m_{i}\cdots  m_{j} \cdots m_{k}\cdots0}\right.\right.\right.}\times\vspace{2mm}\\
&&{\displaystyle\left.\left.\left.\cos \frac{m_{i}\pi
}{l_{i}}x_{i}\cos\frac{m_{j}\pi
}{l_{j}}x_{j}\cos \frac{m_{k}\pi }{l_{k}}x_{k}\right\}\right\}\right\}+\cdots+}\vspace{2mm}\\
&&{\displaystyle\sum^{\infty}_{m_{1}=1}\cdots
\sum^{\infty}_{m_{n}=1}A^{0}_{m_{1}\cdots m_{n}}\cos \frac{m_{1}\pi
}{l_{1}}x_{1}\cdots
\cos \frac{m_{n}\pi }{l_{n}}x_{n}+}\vspace{2mm}\\
&&{\displaystyle\sum^{n}_{i=1}\left\{\sum^{\infty}_{m_{1}=0}\cdots
\sum^{\infty}_{m_{i}=1}
\cdots\sum^{\infty}_{m_{n}=0}A^{1}_{m_{1}\cdots m_{n}}\cos
\frac{m_{1}\pi }{l_{1}}x_{1}\cdots\sin\frac{m_{i}\pi }{l_{i}}x_{i}\cdots \cos \frac{m_{n}\pi }{l_{n}}x_{n}\right\}+}\vspace{2mm}\\
&&{\displaystyle\sum^{n-1}_{i=1}\left\{\sum^{n}_{j=i+1}\left\{\sum^{\infty}_{m_{1}=0}\cdots
\sum^{\infty}_{ m_{i}=1} \cdots\sum^{\infty}_{\stackrel{\ss m_{j}=1}
{i<j}}\cdots\sum^{\infty}_{m_{n}=0}A^{2}_{m_{1}\cdots
m_{n}}\right.\right.}\times\vspace{2mm}\\
&&{\displaystyle\left.\left.\cos \frac{m_{1}\pi
}{l_{1}}x_{1}\cdots\sin\frac{m_{i}\pi }{l_{i}}x_{i}\cdots\sin\frac{m_{j}\pi }{l_{j}}x_{j}\cdots\cos \frac{m_{n}\pi }{l_{n}}x_{n}\right\}\right\}+}\vspace{2mm}\\
&&{\displaystyle\sum^{n-2}_{i=1}\left\{\sum^{n-1}_{j=i+1}\left\{\sum^{n}_{k=j+1}\left\{\sum^{\infty}_{m_{1}=0}\cdots
\sum^{\infty}_{ m_{i}=1} \cdots\sum^{\infty}_{\stackrel{\ss
m_{j}=1}{i<j}} \cdots\sum^{\infty}_{\stackrel{\ss m_{k}=1}
{i<j<k}}\cdots\sum^{\infty}_{m_{n}=0}A^{3}_{m_{1}\cdots
m_{n}}\right.\right.\right.\times}\vspace{2mm}\\
&&{\displaystyle\left.\left.\left.\cos \frac{m_{1}\pi
}{l_{1}}x_{1}\cdots\sin\frac{m_{i}\pi
}{l_{i}}x_{i}\cdots\sin\frac{m_{j}\pi }{l_{j}}x_{j}\cdots
\sin\frac{m_{k}\pi }{l_{k}}x_{k}\cdots\cos \frac{m_{n}\pi }{l_{n}}x_{n}\right\}\right\}\right\}+\cdots+}\vspace{2mm}\\
&&{\displaystyle\sum^{\infty}_{m_{1}=1}\cdots
\sum^{\infty}_{m_{n}=1}A^{n}_{m_{1}\cdots m_{n}}\sin \frac{m_{1}\pi
}{l_{1}}x_{1}\cdots \sin \frac{m_{n}\pi }{l_{n}}x_{n}}.
\vspace{2mm}\\
\end{array}\end{equation}
Noting that, for the fundament solution $K=K(t,x)$, it holds that
\begin{equation}\label{3.8}
\int_{\mathbb{R}^{n}}K(t,x)dx=\int_{\mathbb{R}^{n}}\left(\frac{1}{4\pi
t}\right)^{\frac{n}{2}}\exp\left\{-\frac{|x|^{2}}{4t}\right\}dx=1,
\end{equation}
we obtain
\begin{equation}\label{3.9}\begin{array}{lll}
u(t,x)&=&{\displaystyle\int_{\mathbb{R}^{n}}u_0(\xi)K(t,x-\xi)d\xi}\vspace{2mm}\\
&=&{\displaystyle\int_{\mathbb{R}^{n}}u_0(x-\xi)K(t,\xi)d\xi}\vspace{2mm}\\
&=&{\displaystyle\int_{\mathbb{R}^{n}}u_0(x-\xi)\left(\frac{1}{4\pi
t}\right)^{\frac{n}{2}}\exp\left\{-\frac{|\xi|^{2}}{4t}\right\}d\xi}\vspace{2mm}\\
&=&{\displaystyle\frac{1}{2^{n}l_{1}\cdots
l_{n}}\int^{l_{1}}_{-l_{1}}\cdots\int^{l_{n}}_{-l_{n}}u_0(x)dx_{1}\cdots
dx_{n}+\sum^{n}_{i=1}A^{i}}\vspace{2mm}\\
&=&{\displaystyle\overline{u_0(x)}+\sum^{n}_{i=1}A^{i}},\vspace{2mm}\\
\end{array}\end{equation}
where
\begin{equation}\label{3.9}\begin{array}{lll}
A^{0}&=&{\displaystyle \left(\frac{1}{4\pi
t}\right)^{\frac{n}{2}}\sum^{n}_{i=1}\left\{ \sum^{\infty}_{m_{i}=1}
A^{0}_{0\cdots  m_{i} \cdots 0
}\int_{\mathbb{R}^{n}}\exp\left\{-\frac{\xi_{1}^{2}+\cdots+\xi_{n}^{2}}{4t}\right\}\cos
\frac{m_{i}\pi
}{l_{i}}(x_{i}-\xi_{i})d\xi_{1}\cdots d\xi_{n}\right\}}+\vspace{2mm}\\
&&{\displaystyle\left(\frac{1}{4\pi
t}\right)^{\frac{n}{2}}\sum^{n-1}_{i=1}\left\{\sum^{n}_{j=i+1}\left\{
\sum^{\infty}_{ m_{i}=1} \sum^{\infty}_{\stackrel{\ss m_{j}=1}
{i<j}}A^{0}_{0\cdots  m_{i}\cdots  m_{j} \cdots0
}\int_{\mathbb{R}^{n}}\exp\left\{-\frac{\xi_{1}^{2}+\cdots+\xi_{n}^{2}}{4t}\right\}\times\right.\right.}\vspace{2mm}\\
&&{\displaystyle \left.\left.\cos \frac{m_{i}\pi
}{l_{i}}(x_{i}-\xi_{i})\cos\frac{m_{j}\pi
}{l_{j}}(x_{j}-\xi_{j})d\xi_{1}\cdots d\xi_{n}\right\}\right\}+}\vspace{2mm}\\
&&{\displaystyle \left(\frac{1}{4\pi
t}\right)^{\frac{n}{2}}\sum^{n-2}_{i=1}\left\{\sum^{n-1}_{j=i+1}\left\{\sum^{n}_{k=j+1}\left\{
\sum^{\infty}_{ m_{i}=1} \sum^{\infty}_{\stackrel{\ss m_{j}=1}{i<j}}
\sum^{\infty}_{\stackrel{\ss m_{k}=1} {i<j<k}}A^{0}_{0\cdots  m_{i}\cdots  m_{j} \cdots m_{k}\cdots0}\right.\right.\right.}\times\vspace{2mm}\\
&&{\displaystyle\int_{\mathbb{R}^{n}}\exp\left\{-\frac{\xi_{1}^{2}+\cdots+\xi_{n}^{2}}{4t}\right\}\left.\left.\left.\cos
\frac{m_{i}\pi }{l_{i}}(x_{i}-\xi_{i})\cos\frac{m_{j}\pi
}{l_{j}}(x_{j}-\xi_{j})\cos \frac{m_{k}\pi }{l_{k}}(x_{k}-\xi_{k})\right\}\right\}\right\}+\cdots+}\vspace{2mm}\\
&&{\displaystyle\left(\frac{1}{4\pi
t}\right)^{\frac{n}{2}}\sum^{\infty}_{m_{1}=1}\cdots
\sum^{\infty}_{m_{n}=1}A^{0}_{m_{1}\cdots m_{n}}
\int_{\mathbb{R}^{n}}\exp\left\{-\frac{\xi_{1}^{2}+\cdots+\xi_{n}^{2}}{4t}\right\}\times}\vspace{2mm}\\
&&{\displaystyle
 \cos \frac{m_{1}\pi }{l_{1}}(x_{1}-\xi_{1})\cdots
\cos \frac{m_{n}\pi }{l_{n}}(x_{n}-\xi_{n})d\xi_{1}\cdots d\xi_{n},}\vspace{2mm}\\
\end{array}\end{equation}
\begin{equation}\label{3.9}\begin{array}{lll}
A^{1}&=&{\displaystyle\left(\frac{1}{4\pi
t}\right)^{\frac{n}{2}}\sum^{n}_{i=1}\left\{\sum^{\infty}_{m_{1}=0}\cdots
\sum^{\infty}_{m_{i}=1}
\cdots\sum^{\infty}_{m_{n}=0}A^{1}_{m_{1}\cdots
m_{n}}\int_{\mathbb{R}^{n}}\exp\left\{-\frac{\xi_{1}^{2}+\cdots+\xi_{n}^{2}}{4t}\right\}\times\right.}\vspace{2mm}\\
&&{\displaystyle\left.\cos
\frac{m_{1}\pi }{l_{1}}(x_{1}-\xi_{1})\cdots\sin\frac{m_{i}\pi }{l_{i}}(x_{i}-\xi_{i})\cdots \cos \frac{m_{n}\pi }{l_{n}}(x_{n}-\xi_{n})d\xi_{1}\cdots d\xi_{n}\right\},}\vspace{2mm}\\
\end{array}\end{equation}
\begin{equation}\label{3.9}\begin{array}{lll}
A^{2}&=&{\displaystyle\left(\frac{1}{4\pi
t}\right)^{\frac{n}{2}}\sum^{n-1}_{i=1}\left\{\sum^{n}_{j=i+1}\left\{\sum^{\infty}_{m_{1}=0}\cdots
\sum^{\infty}_{ m_{i}=1} \cdots\sum^{\infty}_{\stackrel{\ss m_{j}=1}
{i<j}}\cdots\sum^{\infty}_{m_{n}=0}A^{2}_{m_{1}\cdots
m_{n}}\int_{\mathbb{R}^{n}}\exp\left\{-\frac{\xi_{1}^{2}+\cdots+\xi_{n}^{2}}{4t}\right\}\right.\right.}\times\vspace{2mm}\\
&&{\displaystyle\left.\left.\cos \frac{m_{1}\pi
}{l_{1}}(x_{1}-\xi_{1})\cdots\sin\frac{m_{i}\pi
}{l_{i}}(x_{i}-\xi_{i})\cdots\sin\frac{m_{j}\pi }
{l_{j}}(x_{j}-\xi_{j})\cdots\cos \frac{m_{n}\pi }{l_{n}}(x_{n}-\xi_{n})d\xi_{1}\cdots d\xi_{n}\right\}\right\},}\vspace{2mm}\\
\end{array}\end{equation}
\begin{equation}\label{3.9}\begin{array}{lll}
A^{3}&=&{\displaystyle\left(\frac{1}{4\pi
t}\right)^{\frac{n}{2}}\sum^{n-2}_{i=1}\left\{\sum^{n-1}_{j=i+1}\left\{\sum^{n}_{k=j+1}\left\{\sum^{\infty}_{m_{1}=0}\cdots
\sum^{\infty}_{ m_{i}=1} \cdots\sum^{\infty}_{\stackrel{\ss
m_{j}=1}{i<j}} \cdots\sum^{\infty}_{\stackrel{\ss m_{k}=1}
{i<j<k}}\cdots\sum^{\infty}_{m_{n}=0}A^{3}_{m_{1}\cdots
m_{n}}\right.\right.\right.\times}\vspace{2mm}\\
&&{\displaystyle\int_{\mathbb{R}^{n}}\exp\left\{-\frac{\xi_{1}^{2}+\cdots+\xi_{n}^{2}}{4t}\right\}\cos
\frac{m_{1}\pi }{l_{1}}(x_{1}-\xi_{1})\cdots\sin\frac{m_{i}\pi
}{l_{i}}(x_{i}-\xi_{i})\cdots}\vspace{2mm}\\
&&{\displaystyle\left.\left.\left.\sin\frac{m_{j}\pi
}{l_{j}}(x_{j}-\xi_{j})\cdots
\sin\frac{m_{k}\pi }{l_{k}}(x_{k}-\xi_{k})\cdots\cos \frac{m_{n}\pi }{l_{n}}(x_{n}-\xi_{n})d\xi_{1}\cdots d\xi_{n}\right\}\right\}\right\}}\vspace{2mm}\\
\end{array}\end{equation}
and so on, particularly,
\begin{equation}\label{3.9}\begin{array}{lll}
A^{n}&=&{\displaystyle\left(\frac{1}{4\pi
t}\right)^{\frac{n}{2}}\sum^{\infty}_{m_{1}=1}\cdots
\sum^{\infty}_{m_{n}=1}A^{n}_{m_{1}\cdots
m_{n}}\int_{\mathbb{R}^{n}}\exp\left\{-\frac{\xi_{1}^{2}+\cdots+\xi_{n}^{2}}{4t}\right\}\times}\vspace{2mm}\\
&&{\displaystyle\sin \frac{m_{1}\pi }{l_{1}}(x_{1}-\xi_{1})\cdots
\sin \frac{m_{n}\pi }{l_{n}}(x_{n}-\xi_{n})d\xi_{1}\cdots
d\xi_{n}}.\vspace{2mm}\\.
\end{array}\end{equation}
On the one hand,
\begin{equation}\label{3.9}\begin{array}{lll}
&&{\displaystyle\!\!\!\!\!\!\!\!\!\!\!\!\!\!\!\!\!\!\!\!\!\left(\frac{1}{4\pi
t}\right)^{\frac{n}{2}}\int_{\mathbb{R}}\exp\left\{-\frac{y^{2}}{4t}\right\}\cos
\frac{n\pi }{l}(x-y)dy}\vspace{2mm}\\
&=&{\displaystyle\left(\frac{1}{4\pi
t}\right)^{\frac{n}{2}}\int_{\mathbb{R}}\exp\left\{-\frac{y^{2}}{4t}\right\}\left[\cos\frac{n\pi
}{l}x\cos\frac{n\pi }{l}y+\sin\frac{n\pi }{l}x\sin\frac{n\pi
}{l}y\right]dy}\vspace{2mm}\\
&=&{\displaystyle\left(\frac{1}{4\pi
t}\right)^{\frac{n}{2}}\cos\frac{n\pi
}{l}x\int_{\mathbb{R}}\exp\left\{-\frac{y^{2}}{4t}\right\}\cos\frac{n\pi
}{l}ydy+\left(\frac{1}{4\pi t}\right)^{\frac{n}{2}}\sin\frac{n\pi
}{l}x\int_{\mathbb{R}}\exp\left\{-\frac{y^{2}}{4t}\right\}\sin\frac{n\pi
}{l}ydy}\vspace{2mm}\\
&=&{\displaystyle\left(\frac{1}{4\pi
t}\right)^{\frac{n}{2}}\cos\frac{n\pi
}{l}x\int_{\mathbb{R}}\exp\left\{-\frac{y^{2}}{4t}\right\}\left[\exp\left\{\frac{n\pi
}{l}yi\right\}+\exp\left\{-\frac{n\pi
}{l}yi\right\}\right]dy+}\vspace{2mm}\\
&&{\displaystyle\frac{1}{2i}\left(\frac{1}{4\pi
t}\right)^{\frac{n}{2}}\sin\frac{n\pi
}{l}x\int_{\mathbb{R}}\exp\left\{-\frac{y^{2}}{4t}\right\}\left[\exp\left\{\frac{n\pi
}{l}yi\right\}-\exp\left\{-\frac{n\pi
}{l}yi\right\}\right]dy}\vspace{2mm}\\
&=&{\displaystyle\left(\frac{1}{4\pi
t}\right)^{\frac{n}{2}}\cos\frac{n\pi
}{l}x\int_{\mathbb{R}}\exp\left\{-\frac{y^{2}}{4t}+\frac{n\pi
}{l}yi\right\}dy+}\vspace{2mm}\\
&&{\displaystyle\left(\frac{1}{4\pi
t}\right)^{\frac{n}{2}}\cos\frac{n\pi
}{l}x\int_{\mathbb{R}}\exp\left\{-\frac{y^{2}}{4t}-\frac{n\pi
}{l}yi\right\}dy+}\vspace{2mm}\\
&&{\displaystyle\left(\frac{1}{4\pi
t}\right)^{\frac{n}{2}}\sin\frac{n\pi
}{l}x\int_{\mathbb{R}}\exp\left\{-\frac{y^{2}}{4t}+\frac{n\pi
}{l}yi\right\}dy-}\vspace{2mm}\\
&&{\displaystyle\frac{1}{2i}\left(\frac{1}{4\pi
t}\right)^{\frac{n}{2}}\sin\frac{n\pi
}{l}x\int_{\mathbb{R}}\exp\left\{-\frac{y^{2}}{4t}-\frac{n\pi
}{l}yi\right\}dy}.\vspace{2mm}\\
\end{array}\end{equation}
On the other hand,
\begin{equation}\label{3.9}\begin{array}{lll}
&&{\displaystyle\!\!\!\!\!\!\!\!\!\!\!\!\!\!\!\!\!\!\!\!\!\!\!\!\!\!\!\!\!\!\!\!\!\!\!\!\!\!\!\!\!\!\int_{\mathbb{R}}\exp\left\{-(ay^{2}+2by+c)\right\}dy}\vspace{2mm}\\
&=&{\displaystyle\int_{\mathbb{R}}\exp\left\{-\frac{1}{a}\left[(ay+b)^{2}+ac-b^2\right]\right\}dy}\vspace{2mm}\\
&=&{\displaystyle
\exp\left\{\frac{b^{2}-ac}{a}\right\}\int_{\mathbb{R}}\exp\left\{-\frac{1}{a}(ay+b)^{2}\right\}dy}\vspace{2mm}\\
&=&{\displaystyle
\left(\frac{4}{a}\right)^{\frac{1}{2}}\exp\left\{\frac{b^{2}-ac}{a}\right\}\int_{0}^{\infty}\exp(-y^{2})dy}\vspace{2mm}\\
&=&{\displaystyle\left(\frac{\pi}{a}\right)^{\frac{1}{2}}\exp\left\{\frac{b^{2}-ac}{a}\right\}\quad (a>0,\;ac-b^{2}>0)}.\vspace{2mm}\\
\end{array}\end{equation}
Thus
\begin{equation}
\left(\frac{1}{4\pi
t}\right)^{\frac{1}{2}}\int_{\mathbb{R}}\exp\left\{-\frac{y^{2}}{4t}\right\}\cos
\frac{n\pi }{l}(x-y)dy=\exp\left\{-\left(\frac{n\pi
}{l}\right)^{2}t\right\}\cos \frac{n\pi }{l}x.
\end{equation}
Similarly,
\begin{equation}
\left(\frac{1}{4\pi
t}\right)^{\frac{1}{2}}\int_{\mathbb{R}}\exp\left\{-\frac{y^{2}}{4t}\right\}\sin
\frac{n\pi }{l}(x-y)dy =\exp\left\{-\left(\frac{n\pi
}{l}\right)^{2}t\right\}\sin \frac{n\pi }{l}x.
\end{equation}
So
\begin{equation}\label{3.9}\begin{array}{lll}
A^{0}&=&{\displaystyle \sum^{n}_{i=1}\left\{ \sum^{\infty}_{m_{i}=1}
A^{0}_{0\cdots  m_{i} \cdots 0 }\left(\frac{1}{4\pi
t}\right)^{\frac{n-1}{2}}\int_{\mathbb{R}^{n-1}}\exp\left\{-\frac{\xi_{1}^{2}+\cdots+\widehat{\xi_{i}^{2}}+\cdots+\xi_{n}^{2}}{4t}\right\}d\xi_{1}\cdots
\widehat{d\xi_{i}}\cdots d\xi_{n}\right.\times}\vspace{2mm}\\
&&{\displaystyle\left. \left(\frac{1}{4\pi
t}\right)^{\frac{1}{2}}\int_{\mathbb{R}}\exp\left\{-\frac{\xi_{i}^{2}}{4t}\right\}\cos
\frac{m_{i}\pi
}{l_{i}}(x_{i}-\xi_{i})d\xi_{i}\right\}}+\vspace{2mm}\\
&&{\displaystyle\sum^{n-1}_{i=1}\left\{\sum^{n}_{j=i+1}\left\{
\sum^{\infty}_{ m_{i}=1} \sum^{\infty}_{\stackrel{\ss m_{j}=1}
{i<j}}A^{0}_{0\cdots  m_{i}\cdots  m_{j} \cdots0
}\left(\frac{1}{4\pi
t}\right)^{\frac{n-2}{2}}\int_{\mathbb{R}^{n-2}}\exp\left\{-\frac{\xi_{1}^{2}+\cdots+\widehat{\xi_{i}^{2}}+\cdots+\widehat{\xi_{j}^{2}}
+\cdots+\xi_{n}^{2}}{4t}\right\}\times\right.\right.}\vspace{2mm}\\
&&{\displaystyle \left.\left.d\xi_{1}\cdots
\widehat{d\xi_{i}}\cdots\widehat{d\xi_{j}}\cdots
d\xi_{n}\left(\frac{1}{4\pi
t}\right)^{\frac{1}{2}}\int_{\mathbb{R}}\exp\left\{-\frac{\xi_{i}^{2}}{4t}\right\}\cos
\frac{m_{i}\pi }{l_{i}}(x_{i}-\xi_{i})d\xi_{i}\right.\right.\times}\vspace{2mm}\\
&&{\displaystyle\left.\left. \left(\frac{1}{4\pi
t}\right)^{\frac{1}{2}}\int_{\mathbb{R}}\exp\left\{-\frac{\xi_{j}^{2}}{4t}\right\}\cos
\frac{m_{j}\pi }{l_{j}}(x_{j}-\xi_{j})d\xi_{j}\right\}\right\}+}\vspace{2mm}\\
&&{\displaystyle
\sum^{n-2}_{i=1}\left\{\sum^{n-1}_{j=i+1}\left\{\sum^{n}_{k=j+1}\left\{
\sum^{\infty}_{ m_{i}=1} \sum^{\infty}_{\stackrel{\ss m_{j}=1}{i<j}}
\sum^{\infty}_{\stackrel{\ss m_{k}=1} {i<j<k}}A^{0}_{0\cdots
m_{i}\cdots  m_{j} \cdots
m_{k}\cdots0}\right.\right.\right.\left(\frac{1}{4\pi
t}\right)^{\frac{n-3}{2}}}\times\vspace{2mm}\\
&&{\displaystyle\int_{\mathbb{R}^{n-3}}\exp\left\{-\frac{\xi_{1}^{2}+\cdots+\widehat{\xi_{i}^{2}}+\cdots+\widehat{\xi_{j}^{2}}
+\cdots+\widehat{\xi_{k}^{2}}+\cdots+\xi_{n}^{2}}{4t}\right\}d\xi_{1}\cdots
\widehat{d\xi_{i}}\cdots\widehat{d\xi_{j}}\cdots\widehat{d\xi_{k}}\cdots
d\xi_{n}\times}\vspace{2mm}\\
&&{\displaystyle\left(\frac{1}{4\pi
t}\right)^{\frac{1}{2}}\int_{\mathbb{R}}\exp\left\{-\frac{\xi_{i}^{2}}{4t}\right\}\cos
\frac{m_{i}\pi }{l_{i}}(x_{i}-\xi_{i})d\xi_{i} \left(\frac{1}{4\pi
t}\right)^{\frac{1}{2}}\int_{\mathbb{R}}\exp\left\{-\frac{\xi_{j}^{2}}{4t}\right\}\cos
\frac{m_{j}\pi }{l_{j}}(x_{j}-\xi_{j})d\xi_{j}\times}\vspace{2mm}\\
&&{\displaystyle\left.\left.\left.\left(\frac{1}{4\pi
t}\right)^{\frac{1}{2}}\int_{\mathbb{R}}\exp\left\{-\frac{\xi_{j}^{2}}{4t}\right\}\cos
\frac{m_{j}\pi }{l_{k}}(x_{k}-\xi_{k})d\xi_{k})\right\}\right\}\right\}+\cdots+}\vspace{2mm}\\
&&{\displaystyle\sum^{\infty}_{m_{1}=1}\cdots
\sum^{\infty}_{m_{n}=1}A^{0}_{m_{1}\cdots m_{n}}\left(\frac{1}{4\pi
t}\right)^{\frac{1}{2}}\int_{\mathbb{R}}\exp\left\{-\frac{\xi_{1}^{2}}{4t}\right\}\cos
\frac{m_{1}\pi
}{l_{1}}(x_{1}-\xi_{1})d\xi_{1}}\times\vspace{2mm}\\
&&{\displaystyle\cdots(\frac{1}{4\pi
t})^{\frac{1}{2}}\int_{\mathbb{R}}\exp\left\{-\frac{\xi_{n}^{2}}{4t}\right\}\cos
\frac{m_{n}\pi }{l_{n}}(x_{n}-\xi_{n})d\xi_{n}}\vspace{2mm}\\
&=&{\displaystyle \sum^{n}_{i=1}\left\{ \sum^{\infty}_{m_{i}=1}
A^{0}_{0\cdots  m_{i} \cdots 0 }\exp\left\{-\left(\frac{m_{i}\pi
}{l_{i}}\right)^{2}t\right\}\cos \frac{m_{i}\pi
}{l_{i}}x_{i}\right\}+}\vspace{2mm}\\
&&{\displaystyle\sum^{n-1}_{i=1}\left\{\sum^{n}_{j=i+1}\left\{
\sum^{\infty}_{ m_{i}=1} \sum^{\infty}_{\stackrel{\ss m_{j}=1}
{i<j}}A^{0}_{0\cdots  m_{i}\cdots  m_{j} \cdots0
}\exp\left\{-\left[\left(\frac{m_{i}\pi
}{l_{i}}\right)^{2}+\left(\frac{m_{j}\pi
}{l_{j}}\right)^{2}\right]t\right\}\cos \frac{m_{i}\pi
}{l_{i}}x_{i}\cos \frac{m_{j}\pi
}{l_{j}}x_{j}\right\}\right\}+}\vspace{2mm}\\
&&{\displaystyle
\sum^{n-2}_{i=1}\left\{\sum^{n-1}_{j=i+1}\left\{\sum^{n}_{k=j+1}\left\{
\sum^{\infty}_{ m_{i}=1} \sum^{\infty}_{\stackrel{\ss m_{j}=1}{i<j}}
\sum^{\infty}_{\stackrel{\ss m_{k}=1} {i<j<k}}A^{0}_{0\cdots
m_{i}\cdots  m_{j} \cdots m_{k}\cdots0}\times
\right.\right.\right.}\vspace{2mm}\\
&&{\displaystyle\left.\left.\left.
\exp\left\{-\left[\left(\frac{m_{i}\pi
}{l_{i}}\right)^{2}+\left(\frac{m_{j}\pi
}{l_{j}}\right)^{2}+\left(\frac{m_{k}\pi
}{l_{k}}\right)^{2}\right]t\right\}\cos \frac{m_{i}\pi
}{l_{i}}x_{i}\cos \frac{m_{j}\pi
}{l_{j}}x_{j}\cos \frac{m_{k}\pi }{l_{k}}x_{k}\right\}\right\}\right\}+\cdots+}\vspace{2mm}\\
&&{\displaystyle\sum^{\infty}_{m_{1}=1}\cdots
\sum^{\infty}_{m_{n}=1}A^{0}_{m_{1}\cdots
m_{n}}\exp\left\{-\left[\left(\frac{m_{1}\pi
}{l_{1}}\right)^{2}+\cdots+\left(\frac{m_{n}\pi
}{l_{n}}\right)^{2}\right]t\right\}\cos
\frac{m_{1}\pi }{l_{1}}x_{1}\cdots \cos \frac{m_{n}\pi }{l_{n}}x_{n},}\vspace{2mm}\\
\end{array}\end{equation}
\begin{equation}\label{3.9}\begin{array}{lll}
A^{1}&=&{\displaystyle\sum^{n}_{i=1}\left\{\sum^{\infty}_{m_{1}=0}\cdots
\sum^{\infty}_{m_{i}=1}
\cdots\sum^{\infty}_{m_{n}=0}A^{1}_{m_{1}\cdots
m_{n}}\left(\frac{1}{4\pi
t}\right)^{\frac{1}{2}}\int_{\mathbb{R}}\exp\left\{-\frac{\xi_{1}^{2}}{4t}\right\}\cos
\frac{m_{1}\pi
}{l_{1}}(x_{1}-\xi_{1})d\xi_{1}\cdots\right.}\vspace{2mm}\\
&&{\displaystyle\left.\left(\frac{1}{4\pi
t}\right)^{\frac{1}{2}}\int_{\mathbb{R}}\exp\left\{-\frac{\xi_{i}^{2}}{4t}\right\}\sin
\frac{m_{i}\pi
}{l_{i}}(x_{i}-\xi_{i})d\xi_{i}\cdots\left(\frac{1}{4\pi
t}\right)^{\frac{1}{2}}\int_{\mathbb{R}}\exp\left\{-\frac{\xi_{n}^{2}}{4t}\right\}\cos
\frac{m_{n}\pi }{l_{n}}(x_{n}-\xi_{n})d\xi_{n}\right\}}\vspace{2mm}\\
&=&{\displaystyle\sum^{n}_{i=1}\left\{\sum^{\infty}_{m_{1}=0}\cdots
\sum^{\infty}_{m_{i}=1}
\cdots\sum^{\infty}_{m_{n}=0}A^{1}_{m_{1}\cdots
m_{n}}\exp\left\{-\left[\left(\frac{m_{1}\pi
}{l_{1}}\right)^{2}+\cdots+\left(\frac{m_{n}\pi
}{l_{n}}\right)^{2}\right]t\right\}\right.}\times\vspace{2mm}\\
&&{\displaystyle\left.\cos \frac{m_{1}\pi }{l_{1}}x_{1}\cdots\sin\frac{m_{i}\pi }{l_{i}}x_{i}\cdots \cos \frac{m_{n}\pi }{l_{n}}x_{n}\right\}},\vspace{2mm}\\
\end{array}\end{equation}
\begin{equation}\label{3.9}\begin{array}{lll}
A^{2}&=&{\displaystyle\sum^{n-1}_{i=1}\left\{\sum^{n}_{j=i+1}\left\{\sum^{\infty}_{m_{1}=0}\cdots
\sum^{\infty}_{ m_{i}=1} \cdots\sum^{\infty}_{\stackrel{\ss m_{j}=1}
{i<j}}\cdots\sum^{\infty}_{m_{n}=0}A^{2}_{m_{1}\cdots
m_{n}}\right.\right.}\times\vspace{2mm}\\
&&{\displaystyle\left(\frac{1}{4\pi
t}\right)^{\frac{1}{2}}\int_{\mathbb{R}}\exp\left\{-\frac{\xi_{1}^{2}}{4t}\right\}\cos
\frac{m_{1}\pi
}{l_{1}}(x_{1}-\xi_{1})d\xi_{1}\cdots\left(\frac{1}{4\pi
t}\right)^{\frac{1}{2}}\int_{\mathbb{R}}\exp\left\{-\frac{\xi_{i}^{2}}{4t}\right\}\sin
\frac{m_{i}\pi }{l_{i}}(x_{i}-\xi_{i})d\xi_{i}\cdots}\vspace{2mm}\\
&&{\displaystyle\left.\left. \left(\frac{1}{4\pi
t}\right)^{\frac{1}{2}}\int_{\mathbb{R}}\exp\left\{-\frac{\xi_{j}^{2}}{4t}\right\}\sin
\frac{m_{j}\pi
}{l_{j}}(x_{j}-\xi_{j})d\xi_{j}\cdots\left(\frac{1}{4\pi
t}\right)^{\frac{1}{2}}\int_{\mathbb{R}}\exp\left\{-\frac{\xi_{n}^{2}}{4t}\right\}\cos
\frac{m_{n}\pi }{l_{n}}(x_{n}-\xi_{n})d\xi_{n}\right\}\right\}}\vspace{2mm}\\
&=&{\displaystyle\sum^{n-1}_{i=1}\left\{\sum^{n}_{j=i+1}\left\{\sum^{\infty}_{m_{1}=0}\cdots
\sum^{\infty}_{ m_{i}=1} \cdots\sum^{\infty}_{\stackrel{\ss m_{j}=1}
{i<j}}\cdots\sum^{\infty}_{m_{n}=0}A^{2}_{m_{1}\cdots
m_{n}}\exp\left\{-\left[\left(\frac{m_{1}\pi
}{l_{1}}\right)^{2}+\cdots+\left(\frac{m_{n}\pi
}{l_{n}}\right)^{2}\right]t\right\}\right.\right.}\times\vspace{2mm}\\
&&{\displaystyle\left.\left.\cos \frac{m_{1}\pi
}{l_{1}}x_{1}\cdots\sin\frac{m_{i}\pi
}{l_{i}}x_{i}\cdots\sin\frac{m_{j}\pi }{l_{j}}x_{j}\cdots\cos
\frac{m_{n}\pi }{l_{n}}x_{n}\right\}\right\}},\vspace{2mm}\\
\end{array}\end{equation}
\begin{equation}\label{3.9}\begin{array}{lll}
A^{3}&=&{\displaystyle\sum^{n-2}_{i=1}\left\{\sum^{n-1}_{j=i+1}\left\{\sum^{n}_{k=j+1}\left\{\sum^{\infty}_{m_{1}=0}\cdots
\sum^{\infty}_{ m_{i}=1} \cdots\sum^{\infty}_{\stackrel{\ss
m_{j}=1}{i<j}} \cdots\sum^{\infty}_{\stackrel{\ss m_{k}=1}
{i<j<k}}\cdots\sum^{\infty}_{m_{n}=0}A^{3}_{m_{1}\cdots
m_{n}}\right.\right.\right.\times}\vspace{2mm}\\
&&{\displaystyle\left(\frac{1}{4\pi
t}\right)^{\frac{1}{2}}\int_{\mathbb{R}}\exp\left\{-\frac{\xi_{1}^{2}}{4t}\right\}\cos
\frac{m_{1}\pi
}{l_{1}}(x_{1}-\xi_{1})d\xi_{1}\cdots\left(\frac{1}{4\pi
t}\right)^{\frac{1}{2}}\int_{\mathbb{R}}\exp\left\{-\frac{\xi_{i}^{2}}{4t}\right\}\sin
\frac{m_{i}\pi }{l_{i}}(x_{i}-\xi_{i})d\xi_{i}\cdots}\vspace{2mm}\\
&&{\displaystyle \left(\frac{1}{4\pi
t}\right)^{\frac{1}{2}}\int_{\mathbb{R}}\exp\left\{-\frac{\xi_{j}^{2}}{4t}\right\}\sin
\frac{m_{j}\pi
}{l_{j}}(x_{j}-\xi_{j})d\xi_{j}\cdots\left(\frac{1}{4\pi
t}\right)^{\frac{1}{2}}\int_{\mathbb{R}}\exp\left\{-\frac{\xi_{i}^{2}}{4t}\right\}\sin
\frac{m_{k}\pi }{l_{k}}(x_{k}-\xi_{k})d\xi_{k}\cdots}\vspace{2mm}\\
&&{\displaystyle \left(\frac{1}{4\pi
t}\right)^{\frac{1}{2}}\int_{\mathbb{R}}\exp\left\{-\frac{\xi_{n}^{2}}{4t}\right\}\cos
\frac{m_{n}\pi }{l_{n}}(x_{n}-\xi_{n})d\xi_{n}}\vspace{2mm}\\
&=&{\displaystyle\sum^{n-2}_{i=1}\left\{\sum^{n-1}_{j=i+1}\left\{\sum^{n}_{k=j+1}\left\{\sum^{\infty}_{m_{1}=0}\cdots
\sum^{\infty}_{ m_{i}=1} \cdots\sum^{\infty}_{\stackrel{\ss
m_{j}=1}{i<j}} \cdots\sum^{\infty}_{\stackrel{\ss m_{k}=1}
{i<j<k}}\cdots\sum^{\infty}_{m_{n}=0}A^{3}_{m_{1}\cdots
m_{n}}\right.\right.\right.\times}\vspace{2mm}\\
&&{\displaystyle \exp\left\{-\left[\left(\frac{m_{1}\pi
}{l_{1}}\right)^{2}+\cdots+\left(\frac{m_{n}\pi
}{l_{n}}\right)^{2}\right]t\right\}}\times\vspace{2mm}\\
&&{\displaystyle\left.\left.\left.\cos \frac{m_{1}\pi
}{l_{1}}x_{1}\cdots\sin\frac{m_{i}\pi
}{l_{i}}x_{i}\cdots\sin\frac{m_{j}\pi
}{l_{j}}x_{j}\cdots\sin\frac{m_{k}\pi }{l_{k}}x_{k}\cdots
\cos \frac{m_{n}\pi }{l_{n}}x_{n}\right\}\right\}\right\}}\vspace{2mm}\\
\end{array}\end{equation}
and so on, in particular,
\begin{equation}\label{3.9}\begin{array}{lll}
A^{n}&=&{\displaystyle\sum^{\infty}_{m_{1}=1}\cdots
\sum^{\infty}_{m_{n}=1}A^{n}_{m_{1}\cdots m_{n}}\left(\frac{1}{4\pi
t}\right)^{\frac{1}{2}}\int_{\mathbb{R}}\exp\left\{-\frac{\xi_{1}^{2}}{4t}\right\}\sin
\frac{m_{1}\pi
}{l_{1}}(x_{1}-\xi_{1})d\xi_{1}\cdots}\vspace{2mm}\\
&&{\displaystyle\left(\frac{1}{4\pi
t}\right)^{\frac{1}{2}}\int_{\mathbb{R}}\exp\left\{-\frac{\xi_{n}^{2}}{4t}\right\}\sin
\frac{m_{n}\pi
}{l_{n}}(x_{n}-\xi_{n})d\xi_{n}}\vspace{2mm}\\
&=&{\displaystyle\sum^{\infty}_{m_{1}=1}\cdots
\sum^{\infty}_{m_{n}=1}A^{n}_{m_{1}\cdots
m_{n}}\exp\left\{-\left[\left(\frac{m_{1}\pi
}{l_{1}}\right)^{2}+\cdots+\left(\frac{m_{n}\pi
}{l_{n}}\right)^{2}\right]t\right\}}\times\vspace{2mm}\\
&&{\displaystyle\sin
\frac{m_{1}\pi }{l_{1}}x_{1}\cdots \sin \frac{m_{n}\pi }{l_{n}}x_{n}}.\vspace{2mm}\\
\end{array}\end{equation}
Noting
\begin{equation}
\exp\left\{-\left[\left(\frac{m_{1}\pi
}{l_{1}}\right)^{2}+\cdots+\left(\frac{m_{n}\pi
}{l_{n}}\right)^{2}\right]t\right\}\longrightarrow 0  \qquad as
\quad t\longrightarrow\infty,
\end{equation}
we obtain the desired (3.2) from (3.9) and (3.19)-(3.24)
immediately. Thus, the proof of Theorem 3.1 is completed.$\quad\Box$

We now prove Theorem 1.2.

\noindent{\bf Proof of Theorem 1.2.} By Theorem 1.1, the Cauchy
problem (\ref{1.2})-(\ref{1.3}) admits a global smooth solution
$x=x(t,\theta_1,\cdots,\theta_n)$, and then by Theorem 3.1, the
solution $x=x(t,\theta)$ satisfies (\ref{1.5}). (\ref{1.5})
implies that hypersurfaces converge to a sphere with radius
$\bar{r}_0$ in the $C^{\infty}$-topology as $t$ goes to the
infinity. Thus, the proof of Theorem 1.2 is completed. $\quad\Box$

\section{Conclusions and open problems}

In this paper we introduce a new geometric flow with rotational
invariance. This flow is described by, formally a system of
hyperbolic partial differential equations with viscosity,
essentially a coupled system of hyperbolic-parabolic partial
differential equations with rotational invariance, which possesses
very interesting geometric properties and dynamical behavior. We
prove that, under this kind of new flow, an arbitrary smooth closed
contractible hypersurface in the Euclidean space
$\mathbb{R}^{n+1}\;(n\ge 1)$ converges to $\mathbb{S}^n$ in the
$C^{\infty}$-topology as $t$ goes to the infinity. As mentioned
before, this result covers the well-known theorem of Gage and
Hamilton in \cite{g-h} for the curvature flow of plane curves and
the famous result of Huisken in \cite{hui} on the flow by mean
curvature of convex surfaces, respectively. In fact, more
applications of this flow to differential geometry and physics can
be expected.

In the present paper, we only investigate the evolution of closed
contractible hypersurfaces in the Euclidean space
$\mathbb{R}^{n}\;(n\ge 2)$ under the flow equation (\ref{1.2}),
there are some fundamental and interesting problems. In particular,
the following open problems seems to us more interesting and
important: (i) use the flow equation (\ref{2.1}) to investigate the
deformation of a closed $m$-dimensional sub-manifold
$x_0=x_0(\theta_1,\cdots,\theta_m)$; (ii) find a suitable way to
extend the results presented in this paper to the case of Riemannian
manifolds in stead of the Euclidean space $\mathbb{R}^{n}$; (iii)
introduce the theory of viscous shock waves to investigate geometric
problems. These problems are worthy to study in the future.

\vskip 0.7cm

 \noindent {\bf\large{Appendix}}
\vskip 0.2cm

\appendix
\renewcommand{\thesection}{A}

In this appendix, we investigate the time-asymptotic behavior of
solution of the following Cauchy problem

\begin{equation}\label{2.5}
\left\{\begin{array}{l}{\displaystyle\frac{\partial u}{\partial
t}+\sum^n_{i=1}\frac{\partial (f_i(u)u)
}{\partial\theta_i}=\Delta u\quad\textrm{in}\,\,(0,\infty)\times {\mathbb{R}}^n},\vspace{2mm}\\
u(0,\theta)= u_0 (\theta) \,\,  \quad\quad\quad\quad\,\,\,\,\quad\textrm{in}\,\, {\mathbb{R}}^n,\\
\end{array}\right.\end{equation}
where $u$ is the unknown scalar function, $f_i(u)\;\;(i=1,\cdots,n)$
are smooth functions, $u_0$ is a periodic function which stands for
the initial data.

\begin{Lemma}
If \,\,$u_0(\theta)$  is a smooth periodic function with period
$\mathbb{T}^{n}=\{\theta=(\theta_{1},\theta_{2},\cdots\cdots
\theta_{n})\mid -l_{i}\leq \theta_{i}\leq l_{i}\}  $  , then
$u(t,x)$ is also periodic with period $\mathbb{T}^{n}$ for any $t$.
\end{Lemma}
\vskip 2mm \noindent{\bf Proof.} Since $u_0(\theta)=u(0,\theta)$ is
periodic with period $\mathbb{T}^{n}$, it holds that
\begin{equation}
u(0,\theta)=u_0(\theta)=u_0(\theta+\mathbb{T}^{n})=u(0,\theta+\mathbb{T}^{n}).
\end{equation}
Taking into account the property of the uniqueness of the solution
$u(t,\theta)$, gives
\begin{equation}
u(t,\theta)=u(t,\theta+\mathbb{T}^{n}).
\end{equation}
This proves Lemma  A.1.$\quad\Box$

Let $\Omega$ be a smooth domain in $\mathbb{R}^{n}$ and consider the
parabolic operator
\begin{equation}
Lu=\sum^n_{i,j=1}a_{ij}(t,x)\frac{\partial^{2}u}{\partial
\theta_{i}\partial
\theta_{j}}+\sum^n_{i=1}b_{i}(t,\theta)\frac{\partial u}{\partial
\theta_{i}}+c(t,\theta)u
\end{equation}
with smooth  and bounded coefficients and a nondegenerate matrix
$(a_{ij})$.
\begin{Lemma} {\bf
(Harnack's inequality)} Suppose that $u(t,\theta)\in
C^{2}((0,T)\times \Omega)$ is a solution of $\frac{\partial
u}{\partial t}-Lu\geq 0$ in $(0,T)\times \Omega$, suppose
furthermore that $u(t,\theta)\geq 0$ in $(0,T)\times \Omega$. Then
for any given compact subset $D$ of $\Omega$ and each $\tau \in
(0,T)$, there exists a positive constant $C_{1}$ depending only on
$D$, $\tau$ and the coefficients of $L$, such that
\begin{equation}
\sup_{K}u(t-\tau,\theta)\leq C_{1} \inf_{K}u(t,\theta).
\end{equation}
\end{Lemma}

\begin{Lemma}
Suppose that $u(t,\theta)$ is a periodic solution of the following
equation
\begin{equation}
\frac{\partial u}{\partial t}+\sum^n_{i=1}\frac{\partial (f_i(u)u)
}{\partial\theta_i}=\Delta u, \quad \forall\;  (t,\theta)\in
(0,\infty)\times{\mathbb{R}}^n
\end{equation}
with a period, say,
$\mathbb{T}^{n}=\{\theta=(\theta_{1},\theta_{2},\cdots\cdots
\theta_{n})\mid -l_{i}\leq \theta_{i}\leq l_{i}\}$, where $f_i\in
L^{\infty}\,\,(i\in N)$. Then there exists a positive constant
$C_{2}$ depending only on $\mathbb{T}^{n}$,$\tau$ and
 $f_i$, such that
\begin{equation}
\|u(t-\tau,\theta)\|_{L^{\infty}( \mathbb{T}^{n})}\leq C_{2}
\|u(t,\theta)\|_{L^{1}( \mathbb{T}^{n})}.
\end{equation}
\end{Lemma}

 \vskip 2mm \noindent{\bf
Proof.} The proof will be divided into two cases.

{\bf Case I: $u(t,\theta)>0$}

 By Harnack inequality, i.e., Lemma
A.2, there exists a positive constant $C_{3}$ depending only on
$\mathbb{T}^{n}$, $\tau$ and
 $f_i$  such that
\begin{equation}
\sup_{\mathbb{T}^{n}}u(t-\tau,\theta)\leq C
\inf_{\mathbb{T}^{n}}u(t,\theta).
\end{equation}
By the mean-value theorem, there exists a point $\theta_{0}\in
\mathbb{T}^{n}$ such that
\begin{equation}
\overline{u(t,\theta)}\stackrel{\triangle}{=}
\frac{1}{|\mathbb{T}^{n}|}\int_{\mathbb{T}^{n}}u(t,\theta)d\theta=
\frac{1}{|\mathbb{T}^{n}|}|\mathbb{T}^{n}|u(t,\theta_{0})\geq
\inf_{\mathbb{T}^{n}} u(t,\theta),
\end{equation}
Therefore, using (A.8), we have
\begin{equation}\label{3.9}\begin{array}{lll}
\frac{1}{|\mathbb{T}^{n}|}\|u(t,\theta)\|_{L^{1}
(\mathbb{T}^{n})}&=&{\displaystyle\overline{u(t,\theta)}\geq \inf_{\mathbb{T}^{n}}u(t,\theta)}\vspace{2mm}\\
&\geq&{\displaystyle\frac{1}{C}\sup_{\mathbb{T}^{n}}u(t-\tau,\theta)\geq
\frac{1}{C}\|u(t-\tau,\theta)\|_{L^{\infty}(\mathbb{T}^{n})}}.
\end{array}\end{equation}
So
\begin{equation}
\|u(t-\tau,\theta)\|_{L^{\infty}(\mathbb{T}^{n})}\leq C_{3}
\|u(t,\theta)\|_{L^{1}(\mathbb{T}^{n})},
\end{equation}
where $C_{3}=\frac{C_{1}}{|\mathbb{T}^{n}|}$.

{\bf Case I\!I:} General case

 Let
\begin{equation}
u(t,\theta)^{+}=\max\left\{u(t,\theta),0\right\},\quad
u(t,\theta)^{-}=-\min\left\{u(t,\theta),0\right\}.
\end{equation}
Then
\begin{equation}
u(t,\theta)=u(t,\theta)^{+}-u(t,\theta)^{-}
\end{equation}
and
\begin{equation}
u(0,\theta)=u(0,\theta)^{+}-u(0,\theta)^{-}.
\end{equation}
By Theorem 2.1, it is easy to see that the solution $u(t,\theta)$
 is unique. Substituting  (A.13) into (A.6) gives
\begin{equation}
\frac{\partial u^{+}}{\partial t}+\sum^n_{i=1}\frac{\partial
(f_i(u)u^{+}) }{\partial\theta_i}-\Delta u^{+}-\left\{\frac{\partial
u^{-}}{\partial t}+\sum^n_{i=1}\frac{\partial (f_i(u)u^{-})
}{\partial\theta_i}+\Delta u^{-}\right\}=0.
\end{equation}

We turn to investigate the  following Cauchy problems
 \begin{equation}\label{2.5}
\left\{\begin{array}{l}{\displaystyle\frac{\partial u^{+}}{\partial
t}+\sum^n_{i=1}\frac{\partial (f_i(u)u^{+})
}{\partial\theta_i}=\Delta u^{+}\quad\textrm{in}\,\,(0,\infty)\times {\mathbb{R}}^n},\vspace{2mm}\\
u^{+}(0,\theta)= u^{+}_0 (\theta) \,\,  \quad \quad\quad\quad\quad\,\,\,\,\quad\textrm{in}\,\, {\mathbb{R}}^n\\
\end{array}\right.\end{equation}
and
\begin{equation}\label{2.5}
\left\{\begin{array}{l}{\displaystyle\frac{\partial u^{-}}{\partial
t}+\sum^n_{i=1}\frac{\partial (f_i(u)u^{-})
}{\partial\theta_i}=\Delta u^{-}\quad\textrm{in}\,\,(0,\infty)\times {\mathbb{R}}^n},\vspace{2mm}\\
u^{-}(0,\theta)= u^{-}_0 (\theta) \,\, \quad
\quad\quad\quad\quad\,\,\,\,\quad\textrm{in}\,\, {\mathbb{R}}^n.
\end{array}\right.\end{equation}

By making use of the method of the proof of  Theorem 2.1, we can
easily prove the Cauchy problems (A.16) and (A.17) admit the unique
non-negative solution, respectively. Noticing (A.11), we have
\begin{equation}
\sup_{\mathbb{T}^{n}}u(t-\tau,\theta)\leq
\sup_{\mathbb{T}^{n}}u^{+}(t-\tau,\theta)\leq
C_{3}\int_{\mathbb{T}^{n}}u^{+}(t,\theta)d\theta
\end{equation}
and
\begin{equation}
\inf_{\mathbb{T}^{n}}u(t-\tau,\theta)\geq
-\sup_{\mathbb{T}^{n}}u^{-}(t-\tau,\theta)\geq
-C_{3}\int_{\mathbb{T}^{n}}u^{-}(t,\theta)d\theta.
\end{equation}
Consequently,
\begin{equation}\label{3.9}\begin{array}{lll}
\|u(t-\tau,\theta)\|_{L^{\infty}(\mathbb{T}^{n})}&\leq&{\displaystyle C_{3}\int_{\mathbb{T}^{n}}(u^{+}(t,\theta)+u^{-}(t,\theta))d\theta}\vspace{2mm}\\
&\leq&{\displaystyle  C_{3}
\int_{\mathbb{T}^{n}}|u(t,\theta)|d\theta\leq
 C_{3}\|u(t,\theta)\|_{L^{1}(\mathbb{T}^{n})}}.
\end{array}\end{equation}
This is the desired estimate (A.7). Thus, the proof of Lemma A.3 is
completed.$\quad\Box$

Let us recall some notions which will be used later. Let
$C^{\infty}_{per}(\mathbb{T}^{n})$ be the space of
$\mathbb{T}^{n}$-periodic functions in $C^{\infty}(\mathbb{R}^{n})$,
and
\begin{equation}\label{3.5}\left\{\begin{array}{lll}
W^{1,\infty}_{per,loc}(\mathbb{R} \times \mathbb{T}^{n}
)&\stackrel{\triangle}{=}&{\displaystyle
\{u(t,\theta)\,\,|\,\,u(t,\theta)\in
W^{1,\infty}_{loc}(\mathbb{R}^{1+n}) \,\,and\,\, u(t,\theta)\,\, is
\,\, \mathbb{T}^{n}-periodic \,\, in \,\, \theta\},}\vspace{2mm}\\
H^{1}_{per}(\mathbb{T}^{n})&\stackrel{\triangle}{=}&{\displaystyle
\left\{u(t,\theta)\,\,|\,\,u(t,\theta)\in\overline{C^{\infty}_{per}(\mathbb{T}^{n})}\bigcap
H^{1}(\mathbb{T}^{n})\right\},\quad
\|\bullet\|_{H^{1}_{per}(\mathbb{T}^{n})}\stackrel{\triangle}{=}\|\bullet\|_{H^{1}(\mathbb{T}^{n})}.} \vspace{2mm}\\
\end{array}\right.\end{equation}
Noticing $(2.10)$\,\,($i.e.,$ $g_{i}(u)=f_{i}(u)u$), we can rewrite
(A.6) as
\begin{equation}
\frac{\partial u}{\partial t}+{\rm div}_{\theta}g(u)=\Delta u.
\end{equation}
In a manner similar to \cite{Dalibard}, we can prove the following
result on stationary solutions of (A.22).
\begin{Proposition}
Let $$g(v(t,\theta))\in W^{1,\infty}_{per,loc}(\mathbb{R}\times
\mathbb{T}^{n} )^{n},\quad {\rm div}_{\theta}g(v(t,\theta))\in
L^{\infty}_{loc}(\mathbb{R}^{1+n}).$$  Suppose that there exist real
numbers $C_{0}>0,m\geq 0$ and $l\in [0,\frac{n+2}{n-2})$ for $n\geq
3$, such that

\begin{equation}
|g^{\prime}_{i}(v(t,\theta))|\leq C_{0}(1+|v|^{m})
\end{equation}
and
\begin{equation}
|{\rm div}_{\theta}g(v(t,\theta))|\leq C_{0}(1+|v|^{l})
\end{equation}
for all $(t,\theta)\in \mathbb{R}\times \mathbb{T}^{n} $. Suppose
furthermore that the couple $(m,l)$ satisfies at least one of the
following conditions
\begin{equation}
m=0,
\end{equation}
\begin{equation}
 \quad l\in [0,1),
\end{equation}
\begin{equation}
 \quad l<\min\left\{\frac{n+2}{n},2\right\} \quad and \quad there \,\, exists \,\, t_{0}\in
\mathbb{R} \,\, such\,\, that \,\, \,\,{\rm
div}_{\theta}g(v(t_{0},\theta))=0 \,\, for\,\, all \,\, \theta\in
T^{n}.
\end{equation}
Then for any fixed $p\in \mathbb{R}$, there exists a unique solution
$v=v(p,\bullet)\in H^{1}_{per}(\mathbb{T}^{n})$ of the problem
\begin{equation}
-\Delta v(p,\theta)+{\rm div}_{\theta}g(v(p,\theta))=0,\quad
\overline{v(p,\bullet)}=p.
\end{equation}
Moveover, $ v(p,\bullet)$ satisfies the growth property: if $p>q$,
then
\begin{equation}
 v(p,\theta)>v(q,\theta),\,\, \forall\,\, \theta\in \mathbb{T}^{n}.
\end{equation}
\end{Proposition}
\begin{Remark}Usually, the problem (A.28) is  called ``cell problem".\end{Remark}

We now state the main result in this section.
\begin{Theorem}
Suppose that $u_{0}(\theta)\in L^{\infty}_{per}(\mathbb{T}^{n})$,
$$u=u(t,\theta)\in C([0,\infty),L^{1}(\mathbb{T}^{n}))\bigcap
L^{\infty}([0,\infty)\times\mathbb{T}^{n})\bigcap
L^{2}_{loc}([0,\infty),H^{1}_{per}(\mathbb{T}^{n}))$$ be the unique
solution of the Cauchy problem
\begin{equation}
\left\{\begin{array}{l} {\displaystyle\frac{\partial u}{\partial
t}+{\rm div}_{\theta}g(u)=\Delta u \quad
\textrm{in}\,\,(0,\infty)\times\mathbb{T}^n },\vspace{2mm}\\
u(\theta,0)= u_0 (\theta) \,\,  \quad
\textrm{in}\,\, \mathbb{T}^n,\\
\end{array}\right.\end{equation}
and $v(\overline{u_0},\theta)\in H^{1}_{per}(\mathbb{T}^{n})$ is the
solution of \,the associated cell problem (A.28), where
$$\overline{u_0}\stackrel{\triangle}{=}\frac{1}{{\rm
vol}\{\mathbb{T}^n\}}\int_{\mathbb{T}^n}u_0(x)dx.$$ Suppose
furthermore  that $g(v(t,\theta))\in
W^{1,\infty}_{per,loc}(\mathbb{R}\times\mathbb{T}^{n} )^{n}$,
$\partial_{\theta_{i}}g^{\prime}_{i}(u(t,\theta))\in
 L^{\infty}_{loc}(R\times \mathbb{T}^{n})$ and
the assumptions of Proposition A.1 are satisfied. Suppose finally
that there exist
 constants $\beta_{1},\beta_{2}\in \mathbb{R}$ such that
\begin{equation}
v(\beta_{1},\theta)\leq u_{0}(\theta)\leq v(\beta_{2},\theta).
\end{equation}
Then it holds that
\begin{equation}
\|u(t,\theta)-v(\overline{u_0},\theta)\|_{L^{\infty}(\mathbb{T}^{n})}\rightarrow
0\quad  as \;\; t\rightarrow \infty.
\end{equation}
\end{Theorem}

Before proving Theorem A.1, we introduce the following notations:
\begin{equation}
\mathbf{M}(t,\theta)\stackrel{\triangle}{=}\sup_{\tau\geq
t}u(\tau,\theta),
\end{equation}
\begin{equation}
\mathbf{P}(t)\stackrel{\triangle}{=}\inf\{p\,\,|\,\,v(p,\theta)\geq
\mathbf{M}(t,\theta)\,\, for\,\, \theta \in \mathbb{T}^{n}\}
\end{equation}
and
\begin{equation}
\mathbf{N}(p,t)\stackrel{\triangle}{=}\{\theta\in
\mathbb{T}^{n}\,\,|\,\,v(p,\theta)< \mathbf{M}(t,\theta)\}.
\end{equation}
By maximum principle, we can find that if $u=u(t,\theta)$ is a
solution of (A.22)  with the initial data satisfying (A.30), then
$$v(\beta_{1},\theta)\leq u(t,\theta)\leq v(\beta_{2},\theta),\quad\forall\;
(t,\theta)\in [0,\infty)\times\mathbb{T}^{n}.
$$
Moreover, it is easy to show that $\mathbf{P}(t)$ is a bounded
non-increasing function of $t$. Therefore, we may set
$$\mathbf{P}\stackrel{\triangle}{=}\lim_{t\rightarrow
\infty}\mathbf{P}(t).$$

 In order to prove Theorem A.1, we need the
following lemma.
\begin{Lemma}
Suppose that $u=u(t,\theta)$ is a solution of
\begin{equation}
\frac{\partial u}{\partial t}+\sum^n_{i=1}\frac{\partial (f_i(u)u)
}{\partial\theta_i}=\Delta u
\end{equation}
and $v(q,\theta)$ is the solution of the associated cell problem
(A.28). Then for any given positive constant $\varepsilon$, there
exist $t_{0}\in \mathbb{R}$ and sequences $\{t_{n}\}$ and
$\{\tau_{n}\}$  satisfying $\tau_{n}\geq t_{n}>t_{0}$ and
${\displaystyle\lim_{n\rightarrow\infty }t_{n}}=\infty $, and
$\theta_{n}\in\mathbf{N}(t_{1},t_{n})$, such that
\begin{equation}
|\omega_{n}(0,\theta_{n} )|\leq \varepsilon,
\end{equation}
where
\begin{equation}
\omega_{n}(t,\theta)\stackrel{\triangle}{=}v(\mathbf{P}(t_{n}),\theta)-u(\tau_{n}+t,\theta)\quad
(t\in[0,1],\theta\in\mathbb{T}^n).
\end{equation}

\end{Lemma}
\vskip 2mm \noindent{\bf Proof.} Since  $v(t,\theta)$ is a
continuous function of $t$, for any given positive constant
$\varepsilon>0$, there exists a positive constant $\delta $ such
that
\begin{equation}
\|v(t,\theta)-v(\mathbf{P},\theta)\|_{L^{\infty}(\mathbb{T}^{n})}\leq
\frac{1}{3}\varepsilon,
\end{equation}
provided that $|t-\mathbf{P}|\leq \delta$. Choose $t_{0}\in
\mathbb{R}$ such that
\begin{equation}
|\mathbf{P}-\mathbf{P}(t)| \leq \delta \quad for \quad  t\geq t_{0}.
\end{equation}
By (A.39), we have
\begin{equation}
\|v(\mathbf{P},\theta)-v(\mathbf{P}(t),\theta)\|_{L^{\infty}(\mathbb
{T}^{n})}\leq \frac{1}{3}\varepsilon.
\end{equation}
Let $t_{1}\in \mathbb{R}$ satisfy
\begin{equation}
|t_{1}-\mathbf{P}| \leq \delta \quad for \quad  t_{1}< \mathbf{P}.
\end{equation}
By (A.39) again, we get
\begin{equation}
\|v(t_{1},\theta)-v(\mathbf{P},\theta)\|_{L^{\infty}(\mathbb{T}^{n})}\leq
\frac{1}{3}\varepsilon.
\end{equation}
Combining (A.41) and (A.43) gives

\begin{equation}\label{3.9}\begin{array}{lll}
\|v(t_{1},\theta)-v(\mathbf{P}(t),\theta)\|_{L^{\infty}(\mathbb{T}^{n})}
&\leq&{\displaystyle
\|v(t_{1},\theta)-v(\mathbf{P},\theta)\|_{L^{\infty}(\mathbb{T}^{n})}+\|v(\mathbf{P},\theta)-v(\mathbf{P}(t),\theta)\|_{L^{\infty}(\mathbb{T}^{n})}}\vspace{2mm}\\
&\leq&{\displaystyle
\frac{1}{3}\varepsilon+\frac{1}{3}\varepsilon=\frac{2}{3}\varepsilon}.\vspace{2mm}\\
\end{array}\end{equation}
Hence, for $t\geq t_{0}$ and $\theta\in \mathbf{N}(t_{1},t)$, it
holds that
\begin{equation}
v(\mathbf{P}(t),\theta)-\frac{2}{3}\varepsilon\leq
v(t_{1},\theta)\leq \mathbf{M}(t,\theta)\leq
v(\mathbf{P}(t),\theta).
\end{equation}
Let $\{t_{n}\}$ be a sequence and satisfy
$$\lim_{n\rightarrow\infty }t_{n}=\infty.$$
For $\theta_{n}\in \mathbf{N}(t_{1},t_{n})$, there exists
$\tau_{n}\,\,(\geq t_{n})$ such that
\begin{equation}
|u(\tau_{n},\theta_{n})-\mathbf{M}(t_{n},\theta_{n})|\leq
\frac{1}{3}\varepsilon.
\end{equation}
Consequently, for large $n$  we have  $t_{n}>t_{0}$, and then we
obtain from (A.45) that
\begin{equation}
v(\mathbf{P}(t_{n}),\theta_{n})-\frac{2}{3}\varepsilon\leq
v(t_{1},\theta_{n})\leq \mathbf{M}(t,\theta_{n})\leq
v(\mathbf{P}(t_{n}),\theta_{n}).
\end{equation}
Combining (A.46) and (A.47) yields
\begin{equation}
\omega_{n}(0,\theta_{n})=|v(\mathbf{P}(t_{n}),\theta_{n})-u(\tau_{n},\theta_{n})|\leq
\varepsilon.
\end{equation}
This proves the desired (A.37). $\quad\Box$

We now prove Theorem A.1 .
 \vskip 2mm \noindent{\bf Proof of Theorem A.1.}
Without loss of generality, we may choose $\tau_{n}$ and $t_{n}$,
with $\tau_{n}\geq t_{n}$, it follows that $\omega_{n}$ is a
non-negative function. By (A.28) and (A.31), $\omega_{n}$ is a
non-negative solution of the following equation
\begin{equation}
\frac{\partial \omega_{n}}{\partial t}+{\rm
div}_{\theta}\left\{\int^{1}_{0}g^{\prime}[\tau
v(\mathbf{P}(t_{n}),\theta)+(1-\tau)u(\tau_{n}+t,\theta)]d\tau
\omega_{n}\right\}=\Delta \omega_{n}\quad
\textrm{in}\,\,(0,\infty)\times\mathbb {T}^n.
\end{equation}
Since $v(t,\theta)$ is a solution of the associated cell problem
(A.28), we can choose $K>0$ such that
\begin{equation}
-K\leq v(t_{1},\theta)\leq v(t_{2},\theta)\leq K.
\end{equation}
So we have
\begin{equation}
\left\|\int^{1}_{0}g^{\prime}[\tau
v(\mathbf{P}(t_{n}),\theta)+(1-\tau)u(\tau_{n}+t,\theta)]d\tau\right\|_{L^{\infty}([0,1]\times\mathbb
{T}^n)}\leq \|g^{\prime}\|_{L^{\infty}([-K,K]\times\mathbb {T}^n)}.
\end{equation}
According to Lemma A.2, there exists a constant $C$ only depending
on $\mathbb {T}^n$ and
$\|\partial_{u}g\|_{L^{\infty}([-K,K]\times\mathbb {T}^n)}$ such
that
\begin{equation}
\sup_{\mathbb {T}^n}\omega_{n}(-\alpha,\theta)\leq C \inf_{\mathbb
{T}^n}\omega_{n}(0,\theta).
\end{equation}
It follows from Lemma A.4 that
\begin{equation}
0\leq v(\mathbf{P}(t_{n}),\theta)-u(\tau_{n}-\alpha,\theta)\leq
C\varepsilon.
\end{equation}
Thus, there exists a sequence $\eta_{n}$ such that
\begin{equation}
\|u(\eta_{n},\theta)-v(\mathbf{P},\theta)\|_{L^{\infty}(\mathbb
{T}^{n})}\rightarrow 0 \quad as \;\; n\rightarrow\infty.
\end{equation}
On the one  hand, integrating equation (A.31) over all $\mathbb
{T}^{n}$ leads to
\begin{equation}
\int_{\mathbb {T}^{n}}u(t,\theta)d\theta=\int_{\mathbb
{T}^{n}}u(0,\theta)d\theta,
\end{equation}
which implies that the total mass of solutions is conserved for all
time. On the other hand, letting $n \rightarrow\infty$ and combining
(A.28) and (A.54)  gives
\begin{equation}
\overline{u_{0}}=\mathbf{P}.
\end{equation}
If $u_{1},u_{2}$ are solutions of (A.36), we can obtain the $L^{1}$
contraction property  by similar method in \cite{Denis},
\begin{equation}
\|u_{1}(t)-u_{2}(t)\|_{L^{1}(\mathbb
{T}^{n})}\leq\|u_{1}(s)-u_{2}(s)\|_{L^{1}(\mathbb {T}^{n})}\quad for
\;\; 0 \leq s\leq t.
\end{equation}
Because $v(\overline{u_{0}},\theta)$ is a stationary solution of
(A.31), we choose $u_{1}=u,u_{2}=v(\overline{u_{0}},\theta),
s=\eta_{n} , t\geq\eta_{n}$, and then we have
\begin{equation}
\|u(t,\theta)-v(\overline{u_{0}},\theta)\|_{L^{1}(\mathbb
{T}^{n})}\rightarrow 0\quad as \;\; t\rightarrow \infty.
\end{equation}
By Lemma A.3 and (A.58), we  obtain
\begin{equation}
\|u(t,\theta)-v(\overline{u_{0}},\theta)\|_{L^{\infty}(\mathbb
{T}^{n})}\rightarrow 0\quad  as \;\; t\rightarrow \infty.
\end{equation}
Thus, the proof of Theorem A.1 is completed. $\quad\Box$

\vskip 10mm

\noindent{\Large {\bf Acknowledgement.}} This work was supported in
part by the NNSF of China (Grant No. 10971190).

\end{document}